\documentclass[11pt,oneside,reqno]{amsart}%
\usepackage{amssymb}
\usepackage{color}
\usepackage{tikz-cd}

\usetikzlibrary{decorations.markings}
\usepackage{xcolor}
\definecolor{deepgreen}{cmyk}{0.99998,0,1,0}

\usepackage{hyperref}
\hypersetup{hypertex=true,
	colorlinks=true,
	linkcolor=blue,
	anchorcolor=blue,
	citecolor=blue}
\usepackage{comment}
\usepackage{mathrsfs}

\newcommand{\cref}{\S\ \ref}
\newcommand{\Cref}{\S\ \ref}

\allowdisplaybreaks[4]

\usepackage{etoolbox}

\usepackage{bm}
\usepackage{bbm}
\theoremstyle{definition}
\newtheorem{defi}{$\mathbf{Definition}$}[section]
\newtheorem*{pro}{$\mathbf{Proof}$}

\theoremstyle{plain}
\newtheorem{theo}[defi]{$\mathbf{Theorem}$}
\newtheorem{lemma}[defi]{$\mathbf{Lemma}$}
\newtheorem{coro}[defi]{$\mathbf{Corollary}$}

\newtheorem{prop}[defi]{$\mathbf{Proposition}$}

\theoremstyle{remark}
\newtheorem{remark}[defi]{$\mathbf{Remark}$}

\newtheorem*{oseledec*}{Oseledec Theorem}

\newcommand{\tro}{\mathrm{Tr}}

\newcommand{\pa}{\partial}

\newcommand{\lk}{\left(}
\newcommand{\rk}{\right)}

\newcommand{\lv}{\left\vert}
\newcommand{\rv}{\right\vert}
\newcommand{\lV}{\left\Vert}
\newcommand{\rV}{\right\Vert}

\newcommand{\bv}{\big\vert}

\newcommand{\bbv}{\bigg\vert}

\newcommand{\bV}{\big\Vert}
\newcommand{\BV}{\Big\Vert}
\newcommand{\bbV}{\bigg\Vert}

\newcommand{\op}{\mathrm{Op}}

\definecolor{pink}{RGB}{249,164,186}
\definecolor{grassgreen}{RGB}{128,255,0}

\numberwithin{equation}{section}

\title{Mixed quantization and partial hyperbolicity}

\author{Snir Ben Ovadia$^{1}$}
\author{Qiaochu Ma$^{2,1}$}
\author{Federico Rodriguez-Hertz$^{1}$}
\address{$^1$Department of Mathematics, The Pennsylvania State University}
\address{$^2$Department of Mathematics, Texas A\&M University}

\date{}

\makeatletter
\newcommand*{\rom}[1]{\expandafter\@slowromancap\romannumeral #1@}
\makeatother

\begin{document} 
	\clearpage\maketitle
	\begin{abstract}
		We establish stable quantum ergodicity for spin Hamiltonians, also known as Pauli-Schrödinger operators. Our approach combines new analytic techniques of mixed quantization, inspired by local index theory, with stable ergodicity results for partially hyperbolic systems.
	\end{abstract}

\section{Introduction}

\subsection{Backgrounds}

The quantum ergodicity theorem (QE) of Shnirelman \cite{MR0402834}, Colin de Verdière \cite{MR818831} and Zelditch \cite{MR916129} asserts that for a compact Riemannian manifold with ergodic geodesic flow, the Laplacian has a density one subsequence of eigenfunctions that tends to be equidistributed.

In \cite{MR2838248, MR3615411}, Bismut-Ma-Zhang obtained the asymptotics of analytic torsions for a series of highly \emph{nonunitary} flat vector bundles. As a global spectral invariant, the analytic torsion is generally difficult to calculate explicitly. To address this, they introduced a seminal technique of differential operators with coefficients in Toeplitz operators \cite[\S\,9]{MR3615411}. This was further developed by Puchol \cite{MR4611826} for asymptotic holomorphic torsions and Ma \cite{MR4665497} for full asymptotic torsions.

Inspired by the framework and the technique in \cite{MR2838248,MR3615411,MR4611826,MR4665497}, Ma-Ma \cite{Ma-Ma} established a uniform version of QE for a series of \emph{unitary} flat vector bundles, where the corresponding dynamical system is the horizontal geodesic flow. Later, a similar but weaker, non-uniform version of QE was also obtained by Cekić-Lefeuvre \cite[Theorem 5.1.7]{cekić2024semiclassicalanalysisprincipalbundles}.

In this article, we present a more general version of QE. Roughly speaking, we show that if a small perturbation is introduced to the setup in \cite{Ma-Ma}, the QE still holds for the perturbed operator, and the corresponding perturbed flow remains ergodic. In other words, the results in \cite{Ma-Ma} are stable both analytically and dynamically. Below, we provide a detailed explanation.

\subsection{QE}\label{s1.2}

First, we introduce our main analytic results on QE.

\subsubsection{Geometric setup}\label{s1.2.1}

Let $(X,g^{TX})$ be an $r$-dimensional compact Riemannian manifold with the Levi-Civita connection $\nabla^{TX}$ and the volume form $dv_X$ induced by $g^{TX}$.

Let $N$ be a compact complex manifold and $(L,h^L)$ a positive line bundle over $N$. Let $g^{TN}$ be the Kähler metric on $TN$ induced by the first Chern form $c_1(L,h^L)$ and $dv_N$ the associated volume form. For $p\in\mathbb{N}^*$, let $L^p=L^{\otimes p}$, the $p$-th tensor power of $L$, and $H^{(0,0)}(N,L^p)$ the space of holomorphic sections of $L^p$ over $N$.

Suppose that $\pi_1(X)$ acts holomorphically on $N$, and this action lifts to a holomorphic isomorphism on $L$. Therefore, $H^{(0,0)}(N,L^p)$ is a finite dimensional complex representation of $\pi_1(X)$, and we can define a series of flat vector bundles $\{F_p\}_{p\in\mathbb{N}^*}$ over $X$ by
\begin{equation}\label{m1.1}
	\begin{split}
			F_p&=\pi_1(X)\backslash\big(\widetilde{X}\times H^{(0,0)}(N,L^p)\big),\\
&=\big\{(\widetilde{x},s)\in\widetilde{X}\times H^{(0,0)}(N,L^p)\big\}/\big((\widetilde{x},s)\sim(\gamma \widetilde{x},\gamma s)\ \text{for\ any\ }\gamma\in\pi_1(X)\big).
	\end{split}
\end{equation}
where $\widetilde{X}$ is the universal covering of $X$. Let $\nabla^{F_p}$ be the natural flat connection and $h^{F_p}$ the Hermitian metric induced by the $L^2$-metric on $H^{(0,0)}(N,L^p)$. We form $\mathscr{N}=\pi_1(X)\backslash\big(\widetilde{X}\times N\big)$, a flat $N$-bundle over $X$, and $\mathscr{L}=\pi_1(X)\backslash(\widetilde{X}\times L)$, a line bundle over $\mathscr{N}$. Note that we use the left $\pi_1(X)$-action on $\widetilde{X}$ to be consistent with the notations commonly used in hyperbolic geometry.

We put natural projections $q\colon\mathscr{N}\to X$ and $\pi\colon T^*X\to X$ and let $q^*T^*X$ be the fibre product of $q$ and $\pi$. Let $S^*X$ be the unit cotangent bundle of $X$ and similarly we denote $q^*S^*X$ the pull-back. Clearly we have $q^*T^*X=\pi_1(X)\backslash(T^*\widetilde{X}\times N)$ and $q^*S^*X=\pi_1(X)\backslash(S^*\widetilde{X}\times N)$. We can equip $q^*T^*X$ with a symplectic structure $\omega_{q^*T^*X}$ and a volume form $dv_{q^*T^*X}$ induced by $\omega_{q^*T^*X}$, see \eqref{m3.8} for details.

In the following, we use $(x,\xi)$ to denote a point of $T^*X$ or $S^*X$ and $(x,\xi,w)$ a point of $q^*T^*X$ or $q^*S^*X$. 

\subsubsection{Hamiltonian operator and Hamiltonian flow}

Let $\{\mathscr{H}^{F_p}\}_{p\in\mathbb{N}^*}$ be a series of self-adjoint Hamiltonian operators, where  $\mathscr{H}^{F_p}$ acts on ${C}^\infty(X,F_p)$. We assume that $\mathscr{H}^{F_p}$ belongs to the class of Toeplitz pseudo-differential operators, and we denote by $\mathscr{H}\in{C}^\infty(q^*T^*X)$ its principal symbol, see \cref{s3.2} for details.

We list all eigenvalues $\{\lambda_{p,i}\}_{i\in\mathbb{N}^*}$ of $\mathscr{H}^{F_p}$ with multiplicity and associated orthonormal eigensections $\{u_{p,i}\}_{i\in\mathbb{N}^*}$, that is,
\begin{equation}\label{1.8.'}
	\mathscr{H}^{F_p}u_{p,i}=\lambda_{p,i}u_{p,i},\ \ \ \ \lV u_{p,i}\rV_{L^2(X,F_{p})}^2=1.
\end{equation}

Let $\psi_t^\mathscr{H}$ denote the Hamiltonian flow of $\mathscr{H}$ on $q^*T^*X$ with respect to $\omega_{q^*T^*X}$. For any $c\in\mathbb{R}$, let $\mathscr{H}^{-1}(c)$ be the corresponding level set of $\mathscr{H}$, on which $\psi_t^\mathscr{H}$ also restricts to a flow. Let $dv_{\mathscr{H}^{-1}(c)}$ be the $\psi_t^\mathscr{H}$-invariant Liouville measure on $\mathscr{H}^{-1}(c)$ induced by $dv_{q^*T^*X}$.


Now we can summarize all the geometric objects in the following diagram
\begin{equation}\label{diag1}
	\begin{tikzcd}[ampersand replacement=\&, column sep=small,row sep=normal,background color=white!20]
		F_p=\pi_1(X)\backslash\big(\widetilde{X}\times H^{(0,0)}(N,L^p)\big)\arrow[d,""swap]\& \mathscr{L}^p=\pi_1(X)\backslash\big(\widetilde{X}\times L^p\big)\arrow[l,"Rq_*"swap]\arrow[d,""swap]\\
		X\& \mathscr{N}=\pi_1(X)\backslash\big(\widetilde{X}\times N\big)\arrow[l,"q"swap]\\
		T^*X \arrow[u,"\pi"]\& q^*T^*X=\pi_1(X)\backslash\big(T^*\widetilde{X}\times N\big)\arrow[u,"\pi"]\arrow[l,"q"swap]\\
		S^*X\arrow[u,hook]\&q^*S^*X=\pi_1(X)\backslash\big(S^*\widetilde{X}\times N\big)\arrow[l,"q"swap]\arrow[u,hook]\\
		\ \& \psi_t^\mathscr{H}\curvearrowright\mathscr{H}^{-1}(c)\arrow[u,dashed,"\cong"]\arrow[uu, bend right=85, hook,"i_c"swap]
	\end{tikzcd}
\end{equation}
where $Rq_*$ denotes the direct image.

As an important special case, if $\mathscr{H}^{F_p}=p^{-2}\Delta^{F_p}$, the nonnegative Laplacian, then we have $\mathscr{H}(x,\xi,w)=\Vert\xi\Vert_{T^*X}^2, \mathscr{H}^{-1}(1)=q^*S^*X$ and $\psi_t^\mathscr{H}$ is the \emph{horizontal geodesic flow} $G_t$ on $q^*T^*X$, which is an extension of the geodesic flow $g_t$ on $T^*X$. Also, if $\mathscr{H}^{F_p}$ is a small perturbation of $p^{-2}\Delta^{F_p}$, then $\mathscr{H}$ is a small perturbation of $\Vert\xi\Vert_{T^*X}^2$, and for generic $c>0$, $\mathscr{H}^{-1}(c)$ is diffeomorphic to $q^*S^*X$ through the map
\begin{equation}\label{1.4..}
	\begin{split}
		(x,\xi,w)\in\mathscr{H}^{-1}(c)\mapstochar\dashrightarrow \big(x,\xi/\lV\xi\rV_{T^*X},w\big)\in q^*S^*X.
	\end{split}
\end{equation}
This explains how $\mathscr{H}^{-1}(c)$ fits into the diagram \eqref{diag1} and the dashed arrow means it is not commutative with other arrows.


\subsubsection{Mixed quantization and QE}

Now we briefly describe the \emph{mixed quantization} procedure. 

For $\mathscr{A}\in{C}^\infty\big(q^*T^*X\big)$, its mixed quantization, $\mathrm{Op}_h(T_{p,\mathscr{A}})$, is the image under the composition of the following two maps
\begin{equation}\label{diag2}
	\begin{split}
{C}^\infty\big(q^*T^*X\big)\xrightarrow{T_{p,\cdot}}{C}^\infty\big(T^*X,\pi^*\mathrm{End}(F_p)\big)\xrightarrow{\mathrm{Op}_h(\cdot)}\mathrm{End}\big(L^2(X,F_p)\big).
	\end{split}
\end{equation}
Here $T_{p,\cdot}$ is the Berezin-Toeplitz quantization along the fiber $N$, regulating the behavior of an infinite number of linear spaces, and $\mathrm{Op}_h(\cdot)$ is the Weyl quantization along the base space $T^*X$, governing high-frequency eigensections. Therefore, combining them enables simultaneous control of the high-frequency eigensections of an infinite number of bundles. 

Indeed, the mixed quantization is mainly restricted to specific symbol classes $\cup_{k\in\mathbb{Z}}S^k(q^*T^*X)\subset{C}^\infty\big(q^*T^*X\big)$, see \eqref{b11'} for details. For any $\mathscr{A}\in S^0(q^*T^*X)$, the averaged integral function over the energy surfaces is defined by
\begin{equation}\label{m3.26}
	c\in\mathbb{R}\longmapsto\frac{1}{\mathrm{Vol}(\mathscr{H}^{-1}(c))}\int_{\mathscr{H}^{-1}(c)}\mathscr{A}dv_{\mathscr{H}^{-1}(c)}.
\end{equation}

Now we state our main QE result, and see Theorem \ref{D9} for an integrated version.
\begin{theo}\label{mt1.1}
Suppose that the Hamiltonian flow $\psi_t^\mathscr{H}$ of $\mathscr{H}$ on the energy surface $(\mathscr{H}^{-1}(c),dv_{\mathscr{H}^{-1}(c)})$ is ergodic over a certain interval $c\in [a,b]$. Then there is a subset $\mathbb{B}\subseteq \{(p,j)\in\mathbb{N}^{*,2}\mid a\leqslant \lambda_{p,j}\leqslant b\}$ with the following density one condition
\begin{equation}\label{m1.6}
	\lim_{p\to\infty}\frac{\bv\{j\in\mathbb{N}^*\mid (p,j)\in \mathbb{B}\}\bv}{\bv\{j\in\mathbb{N}^{*}\mid a\leqslant \lambda_{p,j}\leqslant b\}\bv}=1,
\end{equation}
in which eigensections tend to be equidistributed on $\mathscr{H}^{-1}(c)$ for $c\in[a,b]$. That is, for any $\mathscr{A}\in S^0(q^*T^*X)$ such that the function \eqref{m3.26}
is constant for $c\in[a,b]$, we have
\begin{equation}\label{m1.7}
	\begin{split}
	\lim_{p\to \infty}\sup_{\begin{subarray}{c}
		j\in \mathbb{N}^*\\(p,j)\in\mathbb{B}
\end{subarray}}\bbv\langle \op_{p^{-1}}(T_{\mathscr{A},p})u_{p,j},&u_{p,j}\rangle_{L^2(X,F_p)}\\
&-\frac{1}{\mathrm{Vol}(\mathscr{H}^{-1}([a,b]))}\int_{\mathscr{H}^{-1}([a,b])}\mathscr{A}dv_{\mathscr{H}^{-1}([a,b])}\bbv=0.
	\end{split}
\end{equation}
\end{theo}

We are left with a question, under what condition is the dynamical assumption in Theorem \ref{mt1.1} satisfied? This is the main theme of our dynamical part.

\subsection{Ergodic partially hyperbolic flows}



In this subsection, we specialize \eqref{diag1} to
\begin{equation}\label{1.9}
\big(N,c_1(L,h^L)\big)=\big(\mathbb{CP}^n,\omega_{\mathrm{FS}}\big),
\end{equation}
the complex projective space with the Fubini-Study form. Furthermore, we assume that $\pi_1(X)$ acts on $\mathbb{CP}^n$ through a representation
\begin{equation}\label{m1.16..}
	\rho\colon\pi_1(X)\longrightarrow\mathrm{SU}_{n+1}.
\end{equation}






Now we state our main dynamical result, see Theorem \ref{T6.3} for more details.
\begin{theo}\label{mt1.2}
Suppose that $X$ has Anosov geodesic flow and $\rho$ given in \eqref{m1.16..} has dense image in $\mathrm{SU}_{n+1}$. For any open set $q^*S^*X\subset K\subset q^*(T^*X)$, there exist $\varepsilon>0$ and $k\in\mathbb{N}$ such that if $\lv\mathscr{H}-\Vert\xi\Vert_{T^*X}^2\rv_{C^k(K)}\leqslant \varepsilon$, then $\big(\psi_t^\mathscr{H},dv^{}_{\mathscr{H}^{-1}(c)}\big)$ is ergodic in a certain interval $c\in[a,b]$.
\end{theo}

\begin{remark}
By \eqref{1.4..}, the flow $\psi_t^\mathscr{H}$ acts on $q^*S^*X\cong\mathscr{H}^{-1}(c)$, but it is \emph{not} necessarily an extension of a flow along the base $S^*X$, in other words, $\psi_t^\mathscr{H}$ may not preserve the $\mathbb{CP}^n$-fibre structure. The absence of a $\psi_t^\mathscr{H}$-invariant $\mathbb{CP}^n$-foliation poses a major obstacle. Also, it is crucial that the perturbed flow $\psi_t^\mathscr{H}$ is still Hamiltonian, so we retain a $\psi_t^\mathscr{H}$-invariant symplectic form, which is used in several key points in our proof.
\end{remark}



\subsection{Main result}

By combining Theorems \ref{mt1.1} and \ref{mt1.2}, we obtain their intersection, which constitutes our main result.

Let $O_{\mathbb{CP}^n}(-1)$ be the tautological line bundle over $\mathbb{CP}^n$ give by
\begin{equation}
	O_{\mathbb{CP}^n}(-1)=\big\{(w,v)\in \mathbb{CP}^n\times \mathbb{C}^{n+1}\mid v\in w\big\}.
\end{equation}
Let $O_{\mathbb{CP}^n}(1)$ be the dual bundle of $O_{\mathbb{CP}^n}(-1)$ and $O_{\mathbb{CP}^n}(p)=O_{\mathbb{CP}^n}(1)^{\otimes p}$. It is a classical result that
\begin{equation}
	H^{(0,0)}\big(\mathbb{CP}^n,O_{\mathbb{CP}^n}(p)\big)\cong \mathrm{Sym}^p\mathbb{C}^{n+1},
\end{equation}
where $\mathrm{Sym}^p\mathbb{C}^{n+1}$ denotes the $p$-th symmetric tensor product of $\mathbb{C}^{n+1}$, or equivalently,
\begin{equation}
	\mathrm{Sym}^p\mathbb{C}^{n+1}=\mathbb{C}[z_0,\ldots,z_n]_p,
\end{equation}
the space of homogeneous polynomials of degree $p$ with complex variables $(z_0,\ldots,z_n)$, see for example Huybrechts \cite[Proposition 2.4.1]{MR2093043}.

Now in \eqref{diag1} we take
\begin{equation}\label{1.13}
\big(N,L,H^{(0,0)}(N,L^p)\big)=\big(\mathbb{CP}^n,O_{\mathbb{CP}^n}(1),\mathrm{Sym}^p\mathbb{C}^{n+1}\big),
\end{equation}
and $\pi_1(X)$ acts on $(\mathbb{CP}^n,O_{\mathbb{CP}^n}(1))$ through the representation $\rho$ given in \eqref{m1.16..}. Note that the nonnegative Laplacian $p^{-2}\Delta^{F_p}$ in this case is refered as \emph{Pauli-Schrödinger} spin-$p/2$ operator.

\begin{theo}\label{mt1.4}
Under the setup \eqref{diag1} and \eqref{1.13}, we assume that $\rho$ in \eqref{m1.16..} has dense image in $\mathrm{SU}_{n+1}$ and $\mathscr{H}^{F_p}$ is a small perturbation of $p^{-2}\Delta^{F_p}$. Then on the dynamical side, the Hamiltonian flow $\psi_t^{\mathscr{H}}$ of $\mathscr{H}$ on the energy surface $\mathscr{H}^{-1}(c)$ is ergodic for a certain interval $c\in [a,b]$. On the analytic side, there is a subset $\mathbb{B}\subseteq \{(p,j)\in\mathbb{N}^{*,2}\mid a\leqslant \lambda_{p,j}\leqslant b\}$ with density one condition \eqref{m1.6} and the eigensections satisfy the equidistribution property \eqref{m1.7}.
\end{theo}

\begin{remark}
Note that the components of our main result, Theorem \ref{mt1.4}, originate from four distinct areas. Dynamically, along the base, it involves the stable ergodicity of accessible uniformly hyperbolic flows, and along the fiber, it concerns Lie group action stable ergodicity. Analytically, along the base, equidistribution corresponds to semiclassical quantization, while along the fiber, it relates to geometric quantization. This interplay gives the result special interest and leads one to expect deeper connections among these fields.
\end{remark}

\subsection{An example}\label{s1.5.}

To explain Theorem \ref{mt1.4}, we present an example to which it applies, and see \cref{Sa4.2} for more general applications.

We put $(X,N)=(\Gamma_2\backslash \mathbb{H}^2,\mathbb{CP}^1)$ in \eqref{diag1}, a genus $2$ hyperbolic surface, where $\Gamma_2\backslash \mathbb{H}^2$ is a genus $2$ hyperbolic surface, $\Gamma_2\subset\mathrm{PSL}(2,\mathbb{R})$ and 
\begin{equation}\label{1.15;'}
\Gamma_2\cong\{a_1,b_1,a_2,b_2\mid [a_1,b_1]\cdot[a_2,b_2]=1\}.
\end{equation}
We choose an irrational $\theta\in\mathbb{R}$ and set
\begin{equation}\label{0.17'}
	\rho(a_1)=\rho(b_1)=\begin{pmatrix}
		e^{-i\theta\pi/2}      & 0\\
		0    & e^{i\theta\pi/2}  
	\end{pmatrix},\ \ \rho(a_2)=\rho(b_2)=\begin{pmatrix}
		\cos \frac{\theta\pi}{2}&i\sin\frac{\theta\pi}{2} \\
		i\sin\frac{\theta\pi}{2}&\cos \frac{\theta\pi}{2}    
	\end{pmatrix},
\end{equation}
then $\rho$ extends to a representation $\rho\colon \Gamma_2\to\mathrm{SU}(2)$, and $\rho(\Gamma_2)\subset \mathrm{SU}(2)$ is dense since $\theta$ behaves like an Euler angle.

Consider a vector bundle $\Gamma_2\backslash(\mathbb{H}^2\times \mathfrak{su}_2)$, where $\Gamma_2$ acts on $\mathfrak{su}_2$ through \eqref{1.15;'}, \eqref{0.17'} and the adjoint action of $\mathrm{SU}_2$ on $\mathfrak{su}_2$. Let us take $\Theta_1,\Theta_2\in{C}^\infty(\Gamma_2\backslash \mathbb{H}^2,\Gamma_2\backslash (\mathbb{H}^2\times \mathfrak{su}_2))$, and we can view $\Theta_1,\Theta_2\in{C}^\infty(\Gamma_2\backslash \mathbb{H}^2,\Gamma_2\backslash (\mathbb{H}^2\times\mathrm{End}(\mathrm{Sym}^p\mathbb{C}^{2})))$ through the $\mathfrak{su}_2$ action on $\mathrm{Sym}^p\mathbb{C}^{2}$. Let $\eta$ be a $1$-form on $X$ and we define
\begin{equation}
\Theta^{\Gamma_2\backslash(\mathbb{H}^2\times\mathrm{Sym}^p\mathbb{C}^{2})}=p^{-2}\eta\Theta_1\Theta_2\in \mathscr{C}^\infty\big(\Gamma_2\backslash\mathbb{H}^2,T^*(\Gamma_2\backslash\mathbb{H}^2)\otimes\Gamma_2\backslash \big(\mathbb{H}^2\times\mathrm{End}(\mathrm{Sym}^p\mathbb{C}^{2})\big)\big).
\end{equation}
Let $(x_1,x_2)$ be the upper half plane coordinates of $\Gamma_2\backslash\mathbb{H}^2$, then we can form the following perturbed Hamiltonian operator 
\begin{equation}\label{1.22}
	\begin{split}
		\mathscr{H}^{\Gamma_2\backslash(\mathbb{H}^2\times\mathrm{Sym}^p\mathbb{C}^{2})}=-x_2^2\Big(&\big(p^{-1}\nabla^{\Gamma_2\backslash(\mathbb{H}^2\times\mathrm{Sym}^p\mathbb{C}^{2})}_{\pa_{x_1}}-\Theta^{\Gamma_2\backslash(\mathbb{H}^2\times\mathrm{Sym}^p\mathbb{C}^{2}),*}(\pa_{x_1})\big)\\
		&\cdot\big(p^{-1}\nabla^{\Gamma_2\backslash(\mathbb{H}^2\times\mathrm{Sym}^p\mathbb{C}^{2})}_{\pa_{x_1}}+\Theta^{\Gamma_2\backslash(\mathbb{H}^2\times\mathrm{Sym}^p\mathbb{C}^{2})}(\pa_{x_1})\big)\\
		&+\big(p^{-1}\nabla^{\Gamma_2\backslash(\mathbb{H}^2\times\mathrm{Sym}^p\mathbb{C}^{2})}_{\pa_{x_1}}-\Theta^{\Gamma_2\backslash(\mathbb{H}^2\times\mathrm{Sym}^p\mathbb{C}^{2}),*}(\pa_{x_1})\big)\\
		&\ \  \cdot\big(p^{-1}\nabla^{\Gamma_2\backslash(\mathbb{H}^2\times\mathrm{Sym}^p\mathbb{C}^{2})}_{\pa_{x_1}}+\Theta^{\Gamma_2\backslash(\mathbb{H}^2\times\mathrm{Sym}^p\mathbb{C}^{2})}(\pa_{x_1})\big)\Big)
	\end{split}
\end{equation}
in which $\Theta^{\Gamma_2\backslash(\mathbb{H}^2\times\mathrm{Sym}^p\mathbb{C}^{2})}$ plays the role of a field. So long as $\lv\eta\rv_{C^k(X)}$ is small enough for a suitable $k\in\mathbb{N}$, Theorem \ref{mt1.4} can be applied to $\mathscr{H}^{\Gamma_2\backslash(\mathbb{H}^2\times\mathrm{Sym}^p\mathbb{C}^{2})}$.

\subsection{Comparison with related work}

Schrader-Taylor \cite{MR0995750} and Zelditch \cite{MR1183602} studied the QE for the Laplacian of a series of vector bundles associated with a \emph{principal bundle}, under the dynamical assumptions in \cite[Theorem 9.1]{MR0995750} and \cite[(0.1)]{MR1183602}. However, they did not provide non-trivial examples for which these dynamical assumptions hold, as discussed in \cite[\S\,8]{MR0995750} and \cite[(3.20)]{MR1183602}.

Simultaneously with an early version of this manuscript, Cekić-Lefeuvre \cite[Theorem 5.1.10]{cekić2024semiclassicalanalysisprincipalbundles} obtained the QE for the Laplacian on a series of vector bundles associated with a nearly flat \emph{principal bundle}. In this case, the Hamiltonian flow is an \emph{extension} of the geodesic flow, and the corresponding stable ergodicity is established using the denseness of holonomy similar to the flat case. The stability of denseness is based on a flow version of Burns-Wilkinson \cite[Theorem B]{MR1717580}, however, this seems to lack a detailed proof, as discussed in \cite[\S\,5.1.3]{cekić2024semiclassicalanalysisprincipalbundles}.

For comparison, let us recall the example provided in \eqref{1.22}. The connection $(\nabla^{\Gamma_2\backslash(\mathbb{H}^2\times\mathrm{Sym}^p\mathbb{C}^{2})}+p\Theta^{\Gamma_2\backslash(\mathbb{H}^2\times\mathrm{Sym}^p\mathbb{C}^{2})})$ is \emph{not} principal, requiring the development of the symbol class of Toeplitz pseudo-differential operators. Moreover, the corresponding Hamiltonian flow is \emph{not} an extension of any flow, and this introduces essential additional difficulties compared to the case of an extension of the geodesic flow. 

\subsection{Acknowledgment}

We would like to thank Emilio Corso for helpful discussions, Xiaonan Ma for critical reading on an early version of this manuscript, and Stephanne Nonnenmacher for bringing our attention to this subject. We would like to thank Anatole Katok Center for
Dynamical Systems and Geometry for support. Q. M. would like to thank Sherry Gong, Nigel Higson, Zhizhang Xie, Guoliang Yu for hospitality and was supported by the NSF grants DMS-1952669 and DMS-2247322.

\section{Semiclassical and Geometric quantizations}\label{B}
	
	In this section, we review some properties of Weyl and Berezin-Topeliez quantizations. These two quantization formalisms share many similarities, such as product and trace formulas, and can be compared with each other. For more details, see Zworski \cite{MR2952218} and Ma-Marinescu \cite{MR2339952,MR2393271}.
	\subsection{Weyl quantization}\label{Ba}
	
	Let us choose a finite good covering $\{U_i\}$ of $X$ such that any finite intersection $U_{i_{1}}\cap\cdots\cap U_{i _{k}}$ is differentiably contractible, see for instance Bott-Tu \cite[Theorem 5.1]{MR658304}. Let $(F,\nabla^F)$ be a unitary flat vector bundle on $X$ with a flat Hermitian metric $h^F$, then it can be represented by a unitary representation $\rho\colon\pi_1(X)\to \mathrm{End}(\mathbb{C}^{\dim_{\mathbb{C}} F})$ such that
	\begin{equation}\label{2.1}
	F={\pi_1(X)}\backslash\big(\widetilde{X}\times\mathbb{C}^{\dim_\mathbb{C} F}\big).
	\end{equation}
Moreover, we can trivialize locally $F|_{U_i}\cong \mathbb{C}^{\dim_{\mathbb{C}} F}$ such that the transition maps $\phi_{ij}\in{C}^\infty(U_i\cap U_j,\mathrm{End}(\mathbb{C}^{\dim_{\mathbb{C}} F}))$ are constant unitary matrices. Let $\lV\cdot\rV_F$ denote the norm on $F$ induced by $h^F$ and $\lV\cdot\rV_{\mathrm{End}(F)}$ the operator norm on $\mathrm{End}(F)$.
	
	Let $\{\phi_i\}$ be a partition of unity subordinate to $\{U_i\}$. Then we can define the semiclassical Sobolev norm $\lV\cdot\rV_{H_h^k(X,F)}$, for any $0<h\leqslant1, k\geqslant0$ and $u\in{C}^\infty(X,F)$,
	\begin{equation}
		\lV u\rV_{H_h^k(X,F)}^2=
			\sum_{i}\sum_{\beta\in \mathbb{N}^r,\lv\beta\rv\leqslant k}h^{2\lv\beta\rv}\bV\pa_x^\beta(\phi_i u)\bV_{L^2(U_i,F)}^2,
	\end{equation}
	and through the dual with respect to the $L^2$-metric,
	\begin{equation}
		\lV u\rV_{H_h^{-k}(X,F)}^2=\sup_{\lV u'\rV_{H^{k}_h(X,F)}=1}\lv\langle u,u'\rangle_{L^2(X,F)}\rv^2.
	\end{equation}

	Let $\pi\colon T^*X\to X$ be the natural projection. For $d\in\mathbb{Z}$, we define the symbol class $S_F^d$ as the set of smooth sections $A(x,\xi)\in{C}^\infty(T^*X,\pi^*\mathrm{End}(F))$ such that for any $j\in\mathbb{N}$, the following Kohn-Nirenberg norm is finite
	\begin{equation}\label{b2}
		\lv A\rv^{(d)}_j=\max_{\begin{subarray}{c}\beta_1,\beta_2\in \mathbb{N}^r,
				\lv\beta_1\rv,\lv\beta_2\rv\leqslant j,\\
				U_i\in\{U_i\}
		\end{subarray}}\sup_{x\in U_i,\xi\in T^*_xX}\langle\xi\rangle^{-d+\lv\beta_2\rv}\bV\partial_x^{\beta_1}\partial_{\xi}^{\beta_2}A(x,\xi)\bV_{\mathrm{End}(F_x)}.
	\end{equation}
	We then set $S_F^{-\infty}=\bigcap_{k\in\mathbb{Z}}S_F^d$. 
	
	Let $\{\varphi_i\}$ be a set of nonnegative smooth functions such that $\mathrm{supp}(\varphi_i)\subset U_i$ and $\varphi_i\equiv1$ on an open set containing $\mathrm{supp}(\phi_i)$.
	
\begin{defi}
	For $A\in S_F^d$, its Weyl quantization is defined by \begin{equation}\label{a.2}
	\op_h(A)=\sum_i \varphi_i\op_h(\phi_i A) \varphi_i
\end{equation}
where $0<h\leqslant1$ and $\op_h(\phi_i A)$ is the \emph{Weyl quantization} on $\mathbb{R}^r$, that is,
\begin{equation}\label{a.5}
	\big(\op_h(\phi_i A)u\big)(x)=\frac{1}{(2\pi h)^r}\int_{\mathbb{R}^{r}}\int_{\mathbb{R}^{r}}e^{\frac{\sqrt{-1}}{h}\langle x-x',\xi\rangle}\big(\phi_iA\big)(\textstyle{\frac{x+x'}{2}},\xi)u(x')dx'd\xi.
\end{equation}
	\end{defi}
The integral in \eqref{a.5} is in general divergent, and it should be viewed as a map between Sobolev spaces, see \cite[Proposition E.19]{MR3969938} and \cite[Theorem 4.23]{MR2952218}. 
	
\begin{prop}
For any $d,k\in\mathbb{Z}$, there are $C>0,i\in\mathbb{N}$ such that for any $A\in S^d_F$ and $0<h\leqslant1$, we have
\begin{equation}\label{b11}
	\lV\mathrm{Op}_{h}(A)\rV_{H_h^{d+k}(X,F)\to H^{k}_h(X,F)}\leqslant C\lv A\rv_{i}^{(d)}.
\end{equation}
\end{prop}

	Let $h^\infty\Psi_F^{-\infty}$ be the space of smoothing operators, that is, $B_h\in h^\infty\Psi_F^{-\infty}$ if for any $i\in\mathbb{N}$, there is $C>0$ such that
	\begin{equation}
\lV B_h\rV_{H_h^{-i}(X,F)\to H_h^{i}(X,F)}\leqslant Ch^i.
	\end{equation}
	\begin{defi}
The space $\Psi^d_F$ of pseudo-differential operators is defined by
\begin{equation}
	\Psi^{d}_F=\{\op_h(A_h)+B_h\mid A_h\in S^d_F, B_h\in h^\infty\Psi_F^{-\infty}\}.
\end{equation}
	\end{defi}

	Now we list some useful properties of the Weyl quantization.
	
\begin{lemma}
 For $A\in S_F^0$, we have
\begin{equation}\label{b12}
	\mathrm{Op}_{h}(A)^*=\mathrm{Op}_{h}(A^*),
\end{equation}
where the adjoint is with respect to the $L^2$-metric $\langle\cdot,\cdot\rangle_{L^2(X,F)}$.
\end{lemma}

We have the following product formula \cite[Theorem 9.5]{MR2952218} with an estimate on the remainder term.
\begin{prop}
	The space of pseudo-differential operators is an algebra, in particular, for $A\in S_F^{d_1}$ and $A'\in S_F^{d_2}$, we have
\begin{equation}\label{b7}
\mathrm{Op}_{h}(A)\mathrm{Op}_{h}(A')=\sum_{i=0}^{j}h^i\mathrm{Op}_{h}(A\#_iA')+h^{j+1}R_{h,j+1}(A,A').
\end{equation}
where $A\#_iA'\in S^{d_1+d_2-i}_F$ and $R_{h,j+1}(A,A')\in \Psi^{d_1+d_2-j-1}_F$. Indeed, for any $k\in\mathbb{N}$, there are $C>0,\ell\in\mathbb{N}$ such that
\begin{equation}
		\lv A\#_iA'\rv_{k}^{(d_1+d_2-i)}\leqslant C\lv A\rv_{\ell}^{(d_1)}\lv A'\rv_{\ell}^{(d_2)},
\end{equation}
similarly, for any $k\in\mathbb{Z}$, there are $C>0,\ell\in\mathbb{N}$ such that for any $0<h\leqslant1$,
\begin{equation}
		\lV R_{h,j+1}(A,A')\rV_{H^{k}_h(X,F)\to H^{k+d_1+d_2-j-1}_h(X,F)}\leqslant C\lv A\rv_{\ell}^{(d_1)}\lv A'\rv_{\ell}^{(d_2)}.
\end{equation}
Also,
\begin{equation}\label{b7..}
	A\#_0A=AA',\ \ A\#_1A'=\sum_{j=1}^{r}\frac{1}{2\sqrt{-1}}\big({\pa_{\xi_j}}A\cdot{\pa_{x_j}}A'-{\pa_{x_j}}A\cdot{\pa_{\xi_j}}A'\big).
\end{equation}
\end{prop}

	Let $\omega_{T^*X}$ be the canonical symplectic form on $T^*X$ and $dv_{T^*X}={\omega_{T^*X}^r}/{r!}$ the corresponding volume form on $T^*X$, then locally we have
		\begin{equation}
		\omega_{T^*X}=\sum_{j=1}^{r}d\xi_j\wedge dx_j,\ \  dv_{T^*X}=\prod_{j=1}^{r}d\xi_j\wedge dx_j.
	\end{equation}

\begin{prop}
For any $A\in S_F^{-\infty}$, the operator $\op_h(A)$ is of trace class on $L^2(X,F)$ for any $0<h\leqslant1$, and we have a trace formula
\begin{equation}\label{b.14}
	\frac{(2\pi h)^r}{\dim_\mathbb{C} F}\mathrm{Tr}^{L^2(X,F)}\big[\op_h(A)\big]=\int_{T^*X}\frac{1}{\dim_\mathbb{C} F}\tro^{\pi^*F}[A]dv_{T^*X}.
\end{equation}
Moreover, there is $C>0$ such that for any bounded linear operator $T_h\colon L^2(X,F)\to H^{r+1}_h(X,F)$, it is of trace class on $L^2(X,F)$ and we have
\begin{equation}\label{b.16}
	\frac{(2\pi h)^r}{\dim_{\mathbb{C}} F}\lv\mathrm{Tr}^{L^2(X,F)}[T_h]\rv\leqslant C\lV T_h\rV_{L^2(X,F)\to H^{r+1}_h(X,F)},
\end{equation}
and we note that here $r+1$ is not optimal.
\end{prop}

\subsection{Berezin-Toeplitz quantization}\label{Ca}


Let $N$ be a compact complex manifold with complex dimension $\dim_{\mathbb{C}}N=n$. Let $(L,h^L)$ be a positive holomorphic line bundle on $N$ and $c_1(L,h^L)$ its first Chern form. Then let $g^{TN}$ be the associated Kähler metric on $TN$ and $dv_N=c_1(L,h^L)^n/n!$ the induced volume form on $N$.
	
For $p\in\mathbb{N}^*$, let $L^p=L^{\otimes p}$, the $p$-th tensor power of $L$ and $H^{(0,0)}(N,L^p)$ the space of holomorphic sections of $L^p$ over $N$. Let $\langle\cdot,\cdot\rangle_{H^{(0,0)}(N,L^p)}$ be the $L^2$-product on $H^{(0,0)}(N,L^p)$ induced by $(dv_{N},h^{L})$ and $P_p\colon L^2(N,L^p)\rightarrow H^{(0,0)}(N,L^p)$ the associated orthogonal projection.
	
	\begin{defi}\label{Cb2}
		The Berezin-Toeplitz quantization of $f\in {C}^\infty(N)$ is a sequence of linear operators $\{T_{f,p}\in\mathrm{End}(L^2(N,L^p))\}_{p\in{\mathbb{N}}}$ given by $T_{f,p}=P_pfP_p$, in other words, for any $s,s'\in H^{(0,0)}(N,L^p)$, we have
		\begin{equation}\label{c3}
			\langle T_{f,p}s,s'\rangle_{L^2(N,L^p)}=\int_{N}f(w)\langle s(w),s'(w)\rangle_{h^{L^p}}dv_{N}(w).
		\end{equation}
	\end{defi}

By \eqref{c3}, we have the following obvious result.
\begin{lemma}
For any $f\in C^\infty(N)$, we have
\begin{equation}\label{c4}
	\lV T_{f,p}\rV_{\mathrm{End}(H^{(0,0)}(N,L^p))}\leqslant \lv f\rv_{{C}^0(N)},
\end{equation}
also, the adjoint of the quantization is given by
\begin{equation}\label{c4.}
T_{f,p}^*=T_{\overline{f},p}.
\end{equation}
\end{lemma}

By the Hirzebruch-Riemann-Roch theorem \cite[Theorem 1.4.6]{MR2339952}, we have the following asymptotic dimension formula
\begin{equation}\label{cb5}\dim_\mathbb{C} H^{(0,0)}(N,L^p)=\mathrm{Vol}(N)p^n+O(p^{n-1}).\end{equation}

According to the proof of \cite[Lemma 7.2.4]{MR2339952}, we have the following trace formula with an estimate of the remainder.
	\begin{prop}\label{Cb3}
There is $C>0$ such that for any $f\in{C}^\infty(N)$ and ${p\in{\mathbb{N}^*}}$,
		\begin{equation}\label{cb5'}
			\begin{split}
				\lv\frac{1}{\dim_\mathbb{C} H^{(0,0)}(N,L^p)}\tro^{H^{(0,0)}(N,L^p)}[T_{f,p}]-\int_{N}fd\overline{v}_{N}\rv\leqslant Cp^{-1}\lv f\rv_{{C}^{2}(N)}.
			\end{split}
		\end{equation}
	\end{prop}

	Following \cite[Definition 7.2.1]{MR2339952}, we now define Toeplitz operators.
	\begin{defi}
		A Toeplitz operator is a family of operators  $\big\{T_p\in \mathrm{End}(L^2(N,L^p))\big\}_{p\in{\mathbb{N}^*}}$ such that 
		$T_p=P_pT_pP_p$, and that there exists $\{f_j\in{C}^\infty(N)\}_{j\in \mathbb{N}}$ such that for any $k\in\mathbb{N}$, there is $C>0$ such that
		\begin{equation}\label{cb16}
			\bbV T_p-\sum_{j=0}^k p^{-j}T_{f_j,p}\bbV_{\mathrm{End}(H^{(0,0)}(N,L^p))}\leqslant Cp^{-k-1}.
		\end{equation}
	\end{defi}

From the proof of \cite[Theorem 7.4.1]{MR2339952}, we have the following product formula with an estimate of the remainder.
	\begin{prop}\label{Cb7}
The set of Toeplitz operators is an algebra. In particular, for any $k\in\mathbb{N}$, there is $C>0$ such that for any $f,f'\in{C}^\infty(N)$ and ${p\in{\mathbb{N}^*}}$,
		\begin{equation}\label{cb24}
			\BV T_{f,p}T_{f',p}-\sum_{j=0}^{k}p^{-j}T_{C_j\lk f,f'\rk,p}\BV_{\mathrm{End}(H^{(0,0)}(N,L^p))}\leqslant Cp^{-k-1}\lv f\rv_{{C}^{2k+2}(N)}\cdot\lv f'\rv_{{C}^{2k+2}(N)},
		\end{equation}
		where $C_j(\cdot,\cdot)$ is a smooth bidifferential operator of total degree no more than $2j$. Also,
		\begin{equation}\label{cb23}
			C_0(f,f')=ff',\ \ \ C_1(f,f')-C_1(f',f)=\frac{1}{\sqrt{-1}}\{f,f'\}_N,
		\end{equation}
		where $\{\cdot,\cdot\}_N$ denotes the Poisson bracket of the symplectic form $2\pi c_1(L,h^L)$.
	\end{prop}

\section{Mixed quantization and quantum ergodicity}
	
	In this section, we construct the mixed quantization, which involves semiclassical and geometric quantizations, and use it to deduce QE.

	\subsection{Mixed quantization on $q^*T^*X$}\label{s3.1}
	
	For $d\in\mathbb{Z}$, define the symbol class $S^d(q^*T^*X)$ as the set of smooth functions $\mathscr{A}(x,\xi,w)\in{C}^\infty(q^*T^*X)$ such that for any $j\in\mathbb{N}$, the following norm is finite
	\begin{equation}\label{b11'}
		\lv \mathscr{A}\rv^{(d)}_j=\max_{\begin{subarray}{c}
				\beta_1,\beta_2\in\mathbb{N}^r,\lv\beta_1\rv,\lv\beta_2\rv\leqslant j,\\
				U_i\in\{U_i\}
		\end{subarray}}\sup_{x\in U_i,\xi\in T^*_xX}\langle\xi\rangle^{-d+\lv\beta_2\rv}\bv\partial_x^{\beta_1}\partial_{\xi}^{\beta_2}\mathscr{A}(x,\xi,\cdot)\bv_{{C}^j(q^{-1}(x))},
	\end{equation}
	where $q^{-1}(x)\cong N$.
	
	For $\mathscr{A}\in S^d(q^*T^*X)$ and $p\in\mathbb{N}^*$, let $T_{\mathscr{A},p}$ be its fibrewise Berezin-Toeplitz quantization, that is, $T_{\mathscr{A},p}(x,\xi)=T_{\mathscr{A}(x,\xi,\cdot),p}$. Locally we have for $\beta_1,\beta_2\in \mathbb{N}^r$,
	\begin{equation}\label{m3.3}
		\partial_x^{\beta_1}\partial_{\xi}^{\beta_2}\big(T_{\mathscr{A},p}\big)=T_{\partial_x^{\beta_1}\partial_{\xi}^{\beta_2}\mathscr{A},p},
	\end{equation}
which, together with \eqref{b2}, \eqref{c4} and \eqref{b11'}, implies that $T_{\mathscr{A},p}\in S^d_{F_p}$ and for any $j\in\mathbb{N}$,
	\begin{equation}\label{c11}
		\lv T_{\mathscr{A},p}\rv^{(d)}_j\leqslant \lv \mathscr{A}\rv^{(d)}_j.
	\end{equation}

\begin{defi}
We define the space $\Psi^{d}_{T}$ of Toeplitz operator valued pseudo-differential operators as follows, for a series of pseudo-differential operators $\{\mathscr{T}_p\}_{p\in\mathbb{N}^*}$, in which $\mathscr{T}_p$ acts on ${C}^\infty(X,F_p)$, we say that $\{\mathscr{T}_p\}_{p\in\mathbb{N}^*}\in \Psi^{d}_{T}$ if there exists $\{\mathscr{A}_i\in S^d(q^*T^*X)\}_{i\in\mathbb{N}}$ such that for any $j\in\mathbb{N}$ and $\ell\in\mathbb{Z}$, there is $C>0$ with
\begin{equation}\label{m3.4'}
	\bbV\mathscr{T}_p-\sum_{i=0}^{j}p^{-i}\mathrm{Op}_{p^{-1}}(T_{\mathscr{A}_i,p})\bbV_{H_{p^{-1}}^{\ell+d}(X,F_p)\to H_{p^{-1}}^{\ell}(X,F_p)}\leqslant Cp^{-j-1}.
\end{equation}
Let $\Psi^{-\infty}_T=\bigcap_{d\in\mathbb{Z}}\Psi^{d}_{T}$.
\end{defi}


By \eqref{b7}, \eqref{cb24} and \eqref{m3.3}, we have the following product formula.
\begin{prop}
The set of Toeplitz operator valued pseudo-differential operators is an algebra. In particular, for any $\mathscr{A}\in S^{d_1}(q^*T^*X)$ and $\mathscr{A}'\in S^{d_2}(q^*T^*X)$, we have
	\begin{equation}\label{m3.7}
		\begin{split}
			\mathrm{Op}_{p^{-1}}(T_{\mathscr{A},p})\mathrm{Op}_{p^{-1}}(T_{\mathscr{A}',p})=\sum_{i=0}^{j}p^{-i}\mathrm{Op}_{p^{-1}}(T_{\mathscr{A}*_i\mathscr{A}',p})+p^{-j-1}R_{p,j+1}(\mathscr{A},\mathscr{A}'),
		\end{split}
	\end{equation}
	where $\mathscr{A}*_i\mathscr{A}'\in S^{d_1+d_2}$ and $R_{p,j+1}(\mathscr{A},\mathscr{A}')\in \Psi_T^{d_1+d_2}$. Moreover, for any $k\in\mathbb{N}$, there are $C>0,\ell\in\mathbb{N}$ such that
	\begin{equation}
		\lv A\#_iA'\rv_{k}^{(d_1+d_2)}\leqslant C\lv A\rv_{\ell}^{(d_1)}\lv A'\rv_{\ell}^{(d_2)},
	\end{equation}
and for any $k\in\mathbb{Z}$, there are $C>0,\ell\in\mathbb{N}$ such that for any $0<h\leqslant1$,
	\begin{equation}
		\lV R_{h,j+1}(A,A')\rV_{H^{k}_h(X,F)\to H^{k+d_1+d_2}_h(X,F)}\leqslant C\lv A\rv_{\ell}^{(d_1)}\lv A'\rv_{\ell}^{(d_2)}.
	\end{equation}
Also,
	\begin{equation}\label{m3.8.}
\mathscr{A}*_0\mathscr{A}'=\mathscr{A}\mathscr{A}',\ \ \mathscr{A}*_1\mathscr{A}'=C_1(\mathscr{A},\mathscr{A}')+\mathscr{A}\#_1\mathscr{A}'.
	\end{equation}
\end{prop}

\begin{remark}\label{r3.4}
In \eqref{m3.7}, the order of each term in the expansion and the remainder term remains the same, $(d_1+d_2)$. This is because when expanding the product along the fibre, the order does not decrease. But when one of $\mathscr{A}$ and $\mathscr{A}'$ is fibrewisely constant, what is, $\mathscr{A}\in S^{d_1}(T^*X)$ or $\mathscr{A}'\in S^{d_2}(T^*X)$, we have
$\mathscr{A}*_i\mathscr{A}'\in S^{d_1+d_2-i}$ and $R_{p,j+1}(\mathscr{A},\mathscr{A}')\in \Psi_T^{d_1+d_2-j-1}$ since the fibrewise expansion vanishes.
\end{remark}

Let $\omega_{q^*T^*X}$ be the symplectic form on $q^*T^*X$ and $dv_{q^*T^*X}$ the associated volume form locally given by
\begin{equation}\label{m3.8}
	\omega_{q^*T^*X}=\sum_{j=1}^{r}d\xi_j\wedge dx_j+2\pi c_1(L,h^L),\ \ dv_{q^*T^*X}=\prod_{j=1}^{r}d\xi_j\wedge dx_j\frac{c_1(L,h^L)^n}{n!}.
\end{equation}
Equivalently, the symplectic form $(\omega_{T^*\widetilde{X}}+2\pi c_1(L,h^L))$ on $T^*\widetilde{X}\times N$ is $\pi_1(X)$-invariant, and it descents to $\omega_{q^*T^*X}$. Let $\{\cdot,\cdot\}_{q^*T^*X}$ be the Poisson bracket of $\omega_{q^*T^*X}$. 

\begin{lemma}
For any $\mathscr{A}\in S^{d_1}(q^*T^*X)$ and $\mathscr{A}'\in S^{d_1}(q^*T^*X)$, by \eqref{cb23}, \eqref{m3.7}, \eqref{m3.8.} and \eqref{m3.8}, we have
\begin{equation}\label{m3.10}
	\big[\op_{p^{-1}}(T_{\mathscr{A},p}),\op_{p^{-1}}(T_{\mathscr{A}',p})\big]=\frac{p^{-1}}{\sqrt{-1}}\op_{p^{-1}}\big(T_{\{\mathscr{A},\mathscr{A}'\}_{q^*T^*X},p}\big)+O(p^{-2})_{\Psi^{d_1+d_2}_{T}},
\end{equation}
and in the special case described in Remark \ref{r3.4}, the remainder term can be refined to $O(p^{-2})_{\Psi^{d_1+d_2-2}_{T}}$.
\end{lemma}


By \eqref{b.14}, \eqref{cb5'} and \eqref{c11}, we have the following trace formula.
\begin{prop}
For $\mathscr{A}\in S^{-\infty}(q^*T^*X)$, we have
\begin{equation}\label{m3.9}
	\bigg\vert\frac{(2\pi p^{-1})^r}{\dim_{\mathbb{C}} F_p}\mathrm{Tr}^{L^2(X,F_p)}\big[\op_{p^{-1}}(T_{\mathscr{A},p})\big]-\frac{1}{\mathrm{Vol}(N)}\int_{q^*T^*X}\mathscr{A}dv_{q^*T^*X}\bigg\vert\leqslant Cp^{-1}\lv\mathscr{A}\rv^{(0)}_{2}.
\end{equation}
\end{prop}

\subsection{Local Weyl law}\label{s3.2}

Let us consider a Hamiltonian $\{\mathscr{H}^{F_p}\in\Psi^{2}_{T}\}_{p\in\mathbb{N}^*}$, which is a series of self-adjoint differential operators, and let $\mathscr{H}\in S^2(q^*T^*X)$ be its real principal symbol, which corresponds to the term $\mathscr{A}_0$ in \eqref{m3.4'}, that is,
	\begin{equation}\label{b14}
\mathscr{H}^{F_p}=\op_{p^{-1}}(T_{\mathscr{H},p})+O(p^{-1})_{\Psi^{2}_{T}}.
\end{equation}
We further assume that $\mathscr{H}$ is elliptic in the sense that there exists $C>0$ such that
\begin{equation}\label{m3.14;}
\mathscr{H}(x,\xi,w)\geqslant C\lV\xi\rV_{g^{T^*X}}^2
\end{equation}
outside a compact subset of $q^*T^*X$.

We have the following estimate on the resolvent of $\mathscr{H}^{F_p}$.
\begin{prop}\label{B2}
	For any $\ell\in\mathbb{Z}$, there is $C>0$ such that for any $\lambda\in\mathbb{C}\backslash\mathbb{R}$, 
		\begin{equation}\label{2.12g}
			\lV(\mathscr{H}^{F_p}-\lambda)^{-1}\rV_{H_{p^{-1}}^\ell(X,F_p)\to H^{\ell+1}_{p^{-1}}(X,F_p)}\leqslant C\big((1+\vert\lambda\vert)/\lv\mathrm{Im}\lambda\rv\big)^{\vert\ell\vert+1}.
		\end{equation}
\end{prop}

\begin{pro}
Since \eqref{2.12g} takes a weaker form than that in \cite[Remark 14.6]{MR2952218}, we can provide an elementary proof here. By the self-adjointness assumption on $\mathscr{H}^{F_p}$, we have
\begin{equation}
\mathrm{Im}\big(\big\langle(\mathscr{H}^{F_p}-\lambda)u,u\big\rangle_{L^2(X,F_p)}\big)=-\mathrm{Im}\lambda\lV u\rV_{L^2(X,F_p)}^2,
\end{equation}	
and therefore
\begin{equation}\label{m3.14'}
\bV\big(\mathscr{H}^{F_p}-\lambda\big)^{-1}\bV_{L^2(X,F_p)\to L^2(X,F_p)}\leqslant \lv\mathrm{Im}\lambda\rv^{-1}.
\end{equation}
	By \eqref{b14} and \eqref{m3.14;}, we easily obtain
\begin{equation}
\begin{split}
\lV u\rV_{H_{p^{-1}}^2(X,F_p)}\leqslant C\big(\lV(\mathscr{H}^{F_p}-\lambda)u\rV_{L^2(X,F_p)}+(1+\vert \lambda\vert)\lV u\rV_{L^2(X,F_p)}\big),
\end{split}
\end{equation}	
which, together with \eqref{m3.14'}, implies
\begin{equation}\label{m3.16}
	\begin{split}
		\lV(\mathscr{H}^{F_p}-\lambda)^{-1}\rV_{L^2(X,F_p)\to H_{p^{-1}}^2(X,F_p)}\leqslant C\big(1+(1+\vert\lambda\vert)\lv\mathrm{Im}\lambda\rv^{-1}\big),
	\end{split}
\end{equation}	
which is stronger than \eqref{2.12g} for $\ell=0$. For a vector field $V$ over $X$, we have
\begin{equation}
	\begin{split}
p^{-1}\nabla^{F_p}_V(\mathscr{H}^{F_p}-\lambda)^{-1}=&(\mathscr{H}^{F_p}-\lambda)^{-1}p^{-1}\nabla^{F_p}_V\\
&+(\mathscr{H}^{F_p}-\lambda)^{-1}\big[\mathscr{H}^{F_p},p^{-1}\nabla^{F_p}_V\big](\mathscr{H}^{F_p}-\lambda)^{-1}.
	\end{split}
\end{equation}
Since $p^{-1}\nabla^{F_p}_V$ is a fibrewisely constant differential operator, we can apply \eqref{m3.10} to see that $[\mathscr{H}^{F_p},p^{-1}\nabla^{F_p}_V]\in \Psi_T^{2}$. Iterating from \eqref{m3.16}, we get for $\ell\in\mathbb{N}$,
\begin{equation}\label{m3.18}
\lV(\mathscr{H}^{F_p}-\lambda)^{-1}\rV_{H_{p^{-1}}^\ell(X,F_p)\to H^{\ell+2}_{p^{-1}}(X,F_p)}\leqslant C\big((1+\vert\lambda\vert)/\lv\mathrm{Im}\lambda\rv\big)^{\ell+1},
\end{equation}
which is stronger than \eqref{2.12g}.

Now we consider negative norms. By the self-adjointness of $\mathscr{H}^{F_p}$ and \eqref{m3.16}, we get
\begin{equation}\label{m3.19}
	\begin{split}
\bv\big\langle(\mathscr{H}^{F_p}-\lambda)^{-1}u,u'\big\rangle_{L^2(X,F_p)}\bv&=\bv\big\langle u,(\mathscr{H}^{F_p}-\overline{\lambda})^{-1}u'\big\rangle_{L^2(X,F_p)}\bv\\
&\leqslant \lV u\rV_{H_{p^{-1}}^{-1}(X,F_p)} \lV (\mathscr{H}^{F_p}-\lambda)^{-1}u'\rV_{H_{p^{-1}}^{1}(X,F_p)}\\
&\leqslant C\lV u\rV_{H_{p^{-1}}^{-1}(X,F_p)} \lV u'\rV_{L^2(X,F_p)},
	\end{split}
\end{equation}
where we lost some regularity. From this, we obtain
\begin{equation}
\lV(\mathscr{H}^{F_p}-\lambda)^{-1}\rV_{H_{p^{-1}}^{-1}(X,F_p)\to L^2(X,F_p)}\leqslant C(1+\vert\lambda\vert)/\lv\mathrm{Im}\lambda\rv.
\end{equation}
When $\ell\leqslant-2$, combining \eqref{m3.18} and the argument use in \eqref{m3.19}, we get
\begin{equation}
	\lV(\mathscr{H}^{F_p}-\lambda)^{-1}\rV_{H_{p^{-1}}^\ell(X,F_p)\to H^{\ell+2}_{p^{-1}}(X,F_p)}\leqslant C\big((1+\vert\lambda\vert)/\lv\mathrm{Im}\lambda\rv\big)^{-\ell-1}.
\end{equation}
This completes the proof of \eqref{2.12g}.\qed
\end{pro}

Let $\mathscr{S}(\mathbb{R})$ be the Schwartz function space with a series of norms $\{\lv\cdot\rv_\ell\}_{\ell\in\mathbb{N}}$ given by
\begin{equation}
	\lv\phi\rv_{\ell}=\max_{0\leqslant i,j\leqslant \ell}\sup_{\lambda\in\mathbb{R}}\big\vert \lambda^i\pa_\lambda^j\phi(\lambda)\big\vert.
\end{equation}
Now we give the functional calculus of $\mathscr{H}^{F_p}$.
	\begin{theo}\label{B3.}
		For any $j\in\mathbb{N}$, there are $C>0,\ell\in\mathbb{N}$ such that for any $p\in\mathbb{N}^*$ and $\phi\in\mathscr{S}(\mathbb{R})$, we have $\phi(\mathscr{H}^{F_p})\in\Psi^{-\infty}_{T}$, and
		\begin{equation}\label{b15}
			\lV\phi(\mathscr{H}^{F_p})-\op_{p^{-1}}\big(T_{\phi(\mathscr{H}),p}\big)\rV_{L^2(X,F_p)\to L^2(X,F_p)}\leqslant Cp^{-1}\lv\phi\rv_{\ell}.
		\end{equation}
	\end{theo}
	
	\begin{pro}
We follow the approach in the proof of \cite[Theorem 14.9]{MR2952218}. The assertion $\phi(\mathscr{H}^{F_p})\in\Psi^{-\infty}_{T}$ is deduced from Proposition \ref{B2}, Helffer-Sjöstrand formula \cite[Theorem 11.8]{MR2952218} and Beals theorem \cite[Theorem 9.12]{MR2952218}. Note that for the almost analytic continuation $\widetilde{\phi}$ of $\phi$ and any $i\in\mathbb{N}$, there exist $C>0,\ell\in\mathbb{N}$ such that $\vert\overline{\pa}_\lambda\widetilde{\phi}(\lambda)\vert\leqslant C\vert\mathrm{Im}\lambda\vert^i(1+\vert \lambda\vert^2)^{-i}\lv\phi\rv_{\ell}$, hence the term $(1+\vert\lambda\vert)^{\vert\ell\vert+1}$ in \eqref{2.12g} makes no difference.

For the principal symbol of $\phi(\mathscr{H}^{F_p})$, by \eqref{b11}, \eqref{m3.7} and \eqref{b14}, we get
\begin{equation}
	(\mathscr{H}^{F_p}-\lambda)\op_{p^{-1}}\big(T_{(\mathscr{H}-\lambda)^{-1},p}\big)=\mathrm{Id}+O\big(p^{-1}(1+\vert\lambda\vert)^k\vert\mathrm{Im}\lambda\vert^{-k}\big)_{L^2(X,F_p)\to L^{2}(X,F_p)},
\end{equation}	
then from \eqref{m3.14'} we have
\begin{equation}\label{m3.28''}
	\op_{p^{-1}}\big(T_{(\mathscr{H}-\lambda)^{-1},p}\big)=(\mathscr{H}^{F_p}-\lambda)^{-1}+O\big(p^{-1}(1+\vert\lambda\vert)^k\vert\mathrm{Im}\lambda\vert^{-k-1}\big)_{L^2(X,F_p)\to L^{2}(X,F_p)}.
\end{equation}
Now apply the Helffer-Sjöstrand formula again, we obtain \eqref{b15}.  \qed
\end{pro}

	\begin{theo}\label{B4}
		There are $C>0,i,\ell\in\mathbb{N}$ such that for any $0<h\leqslant1,\phi\in\mathscr{S}(\mathbb{R})$ and $\mathscr{A}\in S^{0}(q^*T^*X)$, we have the following local Weyl law
		\begin{equation}\label{b.18}
			\begin{split}
		\bbv\frac{(2\pi p^{-1})^r}{\dim_{\mathbb{C}} F_p}\tro^{L^2(X,F_p)}\big[\phi(\mathscr{H}^{F_p})\op_{p^{-1}}(T_{\mathscr{A},p})\big]-\frac{1}{\mathrm{Vol}(N)}\int_{q^*T^*X}&\phi(\mathscr{H})\mathscr{A}dv_{T^*X}\bbv\\
		\leqslant &Cp^{-1}\lv\phi\rv_{\ell}\lv\mathscr{A}\rv^{(0)}_i.
			\end{split}
		\end{equation}
	\end{theo}
	
	\begin{pro}
Following the proof of \cite[Theorem 15.3]{MR2952218}, by \eqref{b11}, \eqref{b.16}, \eqref{c4}, \eqref{m3.7}, \eqref{m3.9}, \eqref{b15} and the obvious inequality 
	\begin{equation}\label{b18}
	\lv\frac{1}{\dim_{\mathbb{C}} F_p}\tro^{F_p}[T]\rv\leqslant \lV T\rV_{\mathrm{End}(F_p)}
\end{equation}
where $T\in\mathrm{End}(F_p)$, we get \eqref{b.18}.\qed
	\end{pro}

	\begin{coro}\label{B5}
	For any $a>0$, we have the following Weyl law
		\begin{equation}\label{b.19.}
		\lim_{p\to\infty}\frac{(2\pi p^{-1})^r}{\dim_{\mathbb{C}} F_p}\bv\{j\mid 0\leqslant\lambda_{p,j}\leqslant a\}\bv=\frac{1}{\mathrm{Vol}(N)}\mathrm{Vol}\big(\mathscr{H}^{-1}\big([0,a]\big)\big).
		\end{equation}
	\end{coro}

	\begin{pro}
By setting $\mathscr{A}\equiv1$ in Theorem \ref{B4} and approximating the function $\mathbbm{1}_{[0,a]}$ from above and below by functions in $\mathscr{S}(\mathbb{R})$, we get \eqref{b.19.}.
	\end{pro}

\subsection{Egorov theorem}\label{Bc}

Let $\psi^\mathscr{H}_t$ be the Hamiltonian flow of $\mathscr{H}$ with respect to $\omega_{q^*T^*X}$, that is,
\begin{equation}
	\frac{d}{dt}\big(\psi^\mathscr{H}_t\cdot\mathscr{A}\big)=\{\mathscr{H},\mathscr{A}\}_{q^*T^*X}.
\end{equation}

For $p\in\mathbb{N}^*$ and $t\in\mathbb{R}$, we define the Schrödinger propagator $U^{F_p}_{t}$ of $\mathscr{H}^{F_p}$ by
	\begin{equation}\label{b23}
		U^{F_p}_{t}=e^{-\sqrt{-1}tp\mathscr{H}^{F_p}}.
	\end{equation}
The following Egorov theorem asserts that the operator $U^{F_p}_{t}$ is the quantization of $\psi_t^\mathscr{H}$.
\begin{theo}\label{B9}
		There are $j,k\in\mathbb{N}$ such that for any $T>0$, there is $C>0$ such that for any $p\in\mathbb{N}^*,0\leqslant t\leqslant T$ and $\mathscr{A}\in S^{-\infty}(q^*T^*X)$, we have
		\begin{equation}\label{b23.}
			\BV U^{F_p}_{-t}\mathrm{Op}_{p^{-1}}(T_{\mathscr{A},p})U^{F_p}_{t}-\op_{p^{-1}}(T_{\psi^\mathscr{H}_t\cdot \mathscr{A},p})\BV_{L^2(X,F_p)\to L^2(X,F_p)}\leqslant C\lv \mathscr{A}\rv^{(j)}_{k}.
		\end{equation}
	\end{theo}
	
	\begin{pro}
			Follow the proof of \cite[Theorem 15.2]{MR2952218}, we can deduce \eqref{b23.} from \eqref{m3.10} and \eqref{b14}. \qed
	\end{pro}

\subsection{Quantum ergodicity}\label{sM}

For $\mathscr{A}\in S^0(q^*T^*X)$ and $t_0\in\mathbb{R}$, its time average $\langle\mathscr{A}\rangle_{t_0}\in{C}^\infty(q^*T^*X)$ is defined by
\begin{equation}\label{b25'}
	\langle\mathscr{A}\rangle_{t_0}=\frac{1}{t_0}\int_{0}^{t_0}(\psi^\mathscr{H}_t\cdot \mathscr{A})dt,
\end{equation}
and for $0\leqslant a<b<\infty$, we define its quantum variance $\mathrm{Var}^{F_p}_{a,b}(\mathscr{A})$ by 
\begin{equation}\label{m3.36.}
	\begin{split}
		\mathrm{Var}^{F_p}_{a,b}(\mathscr{A})=\frac{(2\pi p^{-1})^r}{\dim_{\mathbb{C}} F_p}\sum_{ap^{2}\leqslant\lambda_{p,j}\leqslant bp^{2}}&\bbv\langle \op_{p^{-1}}(T_{\mathscr{A},p})u_{p,j},u_{p,j}\rangle_{L^2(X,F_p)}\\
		&-\frac{1}{\mathrm{Vol}(\mathscr{H}^{-1}([a,b]))}\int_{\mathscr{H}^{-1}([a,b])}\mathscr{A}dv_{\mathscr{H}^{-1}([a,b])}\bbv^2,
	\end{split}
\end{equation}
where $\{u_{p,i}\}_{i\in\mathbb{N}^*}$ are the orthonormal eigensections of $\mathscr{H}^{F_p}$ as in \eqref{1.8.'}.

We have the following concentration property for $u_{p,j}$ in \eqref{m3.36.}.
\begin{lemma}
There are $k\in\mathbb{N},C>0$ and a compact set $\mathscr{H}^{-1}([a,b])\subset K\subset q^*T^*X$ such that for any $\mathscr{A}\in S^0(q^*T^*X)$ with $\mathrm{supp}(\mathscr{A})\cap K=\emptyset$, we have
\begin{equation}\label{m3.37'}
	\lV\op_{p^{-1}}(T_{\mathscr{A},p})u_{p,j}\rV_{L^2(X,F_p)}\leqslant Cp^{-1}\lv\mathscr{A}\rv^{(0)}_k
\end{equation}
\end{lemma}

\begin{pro}
In \eqref{b14}, the remainder term is in $\Psi_T^2$, so we \emph{cannot} straightforwardly obtain from \eqref{1.8.'} that
\begin{equation}
\Vert\op_{p^{-1}}(T_{(\mathscr{H}-p^{-2}\lambda_{p,j}),p})u_{p,j}\Vert_{L^2(X,F_p)}\leqslant Cp^{-1}.
\end{equation}
Instead, by \eqref{m3.7} we get
	\begin{equation}
		\begin{split}
\op_{p^{-1}}\big(T_{(1+\lv\xi\rv_{g^{T*X}}^2)^{-1}(\mathscr{H}-p^{-2}\lambda_{p,j}),p}\big)=&\op_{p^{-1}}\big(T_{(1+\lv\xi\rv_{g^{T*X}}^2)^{-1},p}\big)\op_{p^{-1}}\big(T_{(\mathscr{H}-p^{-2}\lambda_{p,j}),p}\big)\\
&+O(p^{-1})_{\Psi^0_{T}},
		\end{split}
	\end{equation}
 by \eqref{b11} we have
\begin{equation}\label{m3.39'}
\BV\op_{p^{-1}}\big(T_{(1+\lv\xi\rv_{g^{T*X}}^2)^{-1}(\mathscr{H}-p^{-2}\lambda_{p,j}),p}\big)u_{p,j}\BV_{L^2(X,F_p)}\leqslant Cp^{-1}.
	\end{equation}
Similarly, by \eqref{m3.7} we have
\begin{equation}
	\begin{split}
\op_{p^{-1}}(T_{\mathscr{A},p})=&\op_{p^{-1}}\big(T_{\mathscr{A}(1+\lv\xi\rv_{g^{T*X}}^2)(\mathscr{H}-p^{-2}\lambda_{p,j})^{-1}}\big)	\op_{p^{-1}}\big(T_{(1+\lv\xi\rv_{g^{T*X}}^2)^{-1}(\mathscr{H}-p^{-2}\lambda_{p,j})}\big)\\
&+O(p^{-1}\lv\mathscr{A}\rv^{(0)}_k)_{\Psi^0_{T}},
	\end{split}
\end{equation}
which, together with \eqref{b11} and \eqref{m3.39'}, implies \eqref{m3.37'}.\qed
\end{pro}

Now we give an estimate for $\mathrm{Var}^{F_p}_{a,b}(\mathscr{A})$.
\begin{theo}\label{B10}
	There is $k\in\mathbb{N}$ such that for any $0\leqslant a<b<\infty$ and $0<t_0<\infty$, there exist $C_{a,b}, C_{a,b,t_0}>0$ such that for any $\mathscr{A}\in S^0(q^*T^*X)$ and $p\in\mathbb{N}^*$, we have
	\begin{equation}\label{b26}
		\begin{split}
			\mathrm{Var}^{F_p}_{a,b}(\mathscr{A})\leqslant& C_{a,b}\bbV\langle\mathscr{A}\rangle_{t_0}-\frac{1}{\mathrm{Vol}(\mathscr{H}^{-1}([a,b]))}\int_{\mathscr{H}^{-1}([a,b])}\mathscr{A}dv_{\mathscr{H}^{-1}([a,b])}\bbV^2_{L^2(\mathscr{H}^{-1}([a,b]))}\\
			&+C_{a,b,T}p^{-1}\big(\lv\mathscr{A}\rv^{(0)}_k\big)^2.
		\end{split}
	\end{equation}
\end{theo}

\begin{pro}
By \eqref{m3.37'}, using a cut-off function on $q^*T^*X$, we can suppose without loss of generality that $\mathscr{A}\in S^{-\infty}(q^*T^*X)$. From \eqref{b23}, we have $U^{F_p}_{t}u_{p,j}=e^{-itp^{-1}\lambda_{p,j}}u_{p,j}$, therefore,
	\begin{equation}\label{b27}
		\big\langle U^{F_p}_{-t}\op_{p^{-1}}(T_{\mathscr{A},p})U^{F_p}_{t}u_{p,j},u_{p,j}\big\rangle_{L^2(X,F_p)}=\big\langle\mathrm{Op}_{p^{-1}}(\mathscr{A})u_{p,j},u_{p,j}\big\rangle_{L^2(X,F_p)}.
	\end{equation}
	In \eqref{b.18}, we take a nonnegative $\phi\in\mathscr{S}(\mathbb{R})$ with $\mathrm{supp}(\phi)\subset(0,\infty)$ and  $\phi\equiv1$ on $[a,b]$, then \eqref{b26} follows immediately from \eqref{b11}, \eqref{b12}, \eqref{b7}, \eqref{b.16}, \eqref{b23.} and \eqref{b27}.\qed
\end{pro}

	We have the following integrated form of QE.
	\begin{theo}\label{D9}
			Suppose that the Hamiltonian flow $\psi^\mathscr{H}_t$ is ergodic on $(\mathscr{H}^{-1}(c),d\overline{v}_{\mathscr{H}^{-1}(c)})$ for any $c\in [a,b]$, then for any $\mathscr{A}\in S^0(q^*T^*X)$ satisfies the assumption \eqref{m3.26}, we have
		\begin{equation}\label{m3.27}
			\begin{split}
				\lim_{p\to\infty}\mathrm{Var}^{F_p}_{a,b}(\mathscr{A})=0.
			\end{split}
		\end{equation}
	\end{theo}
	\begin{pro}
Let $p\to \infty$ and then $t_0\to \infty$ in \eqref{b26}, then we obtain \eqref{m3.27} from the mean ergodic theorem. \qed
	\end{pro}
	
Using a diagonal argument to Theorem \ref{D9}, we deduce Theorem \ref{mt1.1}.

\section{Stable ergodicity}

In this section, we present our main dynamical theorem on the ergodicity of a small perturbation of the horizontal geodesic flow.

Throughout this section, we shall work under the setup \eqref{diag1}, \eqref{1.9} and \eqref{m1.16..}, also, we suppose that $\rho$ given in \eqref{m1.16..} has \emph{dense image} and $X$ has \emph{Anosov geodesic flow}. For ease of notations, we use $y=(x,\xi)$ to denote a point of $S^*X$ and $z=(x,\xi,w)$ a point of $q^*S^*X$.

\subsection{Unperturbed dynamical system}\label{S4.1'}


Recall that $g_t\colon S^*X\to S^*X$ denotes the geodesic flow, hence an Anosov flow preserving the Liouville measure $dv_{S^*X}$. Meanwhile, $G_t\colon q^*S^*X\to q^*S^*X$ denotes the horizontal geodesic flow, a $\mathbb{CP}^n$-extension of the geodesic flow $g_t$ that preserves the augmented Liouville measure $dv_{q^*S^*X}$.

We have the following $G_t$ invariant decomposition of the tangent bundle of $q^*S^*X$
\begin{equation}\label{m4.1.}
Tq^*S^*X=E^{c,G_t}\oplus E^{s,G_t}\oplus E^{u,G_t},
\end{equation}
and let $W^{s/u,G_t}_\mathrm{loc}$ be the corresponding local stable and unstable foliations. 

Let us define a vector field $\mathcal{X}_{q^*S^*X}$ over $q^*S^*X$ by the direction of $G_t$, namely
\begin{equation}\label{m4.2.}
	\mathcal{X}_{q^*S^*X}=\frac{d}{dt}\Big|_{t=0}G_t,
\end{equation}
then the central vector bundle $E^{c,G_t}$ further has the following $G_t$-invariant decomposition
\begin{equation}\label{m4.3.}
E^{c,G_t}=\mathbb{C}\cdot\mathcal{X}_{q^*S^*X}\oplus\ker q_*,
\end{equation}
and $\ker q_*$ corresponds to the foliation $\mathcal{F}$ of the $\mathbb{CP}^n$-fibre structure, that is,
\begin{equation}\label{m4.4.}
\mathcal{F}(z)=q^{-1}(q(z))\cong \mathbb{CP}^n.
\end{equation}
The authentic central foliation $\mathcal{F}^{G_t}$ is defined by
\begin{equation}\label{m4.5..}
	\mathcal{F}^{G_t}(z)=\bigcup_{t\in[-\tau,\tau]} G_t\big(\mathcal{F}(z)\big)\cong [-\tau,\tau]\times \mathbb{CP}^n
\end{equation}
where $\tau>0$ is small.

Using the horizontal parallel transport, for $z'\in W^{s/u,G_t}(z)$, we can define the holonomy map $H^{s/u,G_t}_{zz'}\colon \mathcal{F}(z)\to \mathcal{F}(z')$. For a su-path $\gamma^{G_t}=\big(z_1,\cdots,z_{\vert\gamma^{G_t}\vert}\big)$, let us define $R_{\gamma^{G_t}}\colon\mathcal{F}(z_1)\to \mathcal{F}(z_{\vert\gamma^{G_t}\vert})$ by the composition of holonomies
\begin{equation}\label{m4.6.}
R_{\gamma^{G_t}}=H^{s/u,G_t}_{z_{\vert\gamma^{G_t}\vert-1}z_{\vert\gamma^{G_t}\vert}}\cdots H^{s/u,G_t}_{z_1z_2}.
\end{equation}
In particular, if $\gamma^{G_t}$ projects to a su-loop on $S^*X$, in other words, $q(z_1)=q(z_{\vert\gamma^{G_t}\vert})$, then $\mathcal{F}(z_1)=\mathcal{F}(z_{\vert\gamma^{G_t}\vert})$ and we can view $R_{\gamma^{G_t}}\in\mathrm{SU}_{n+1}$ by $\rho$ given in \eqref{m1.16..}. From the dense image assumption on $\rho$, we easily have the following result.
\begin{lemma}
For any $z\in q^*S^*X$, there exist su-paths $\gamma_1^{G_t},\ldots,\gamma_\ell^{G_t}$ from $\mathcal{F}^{}(z)$ to itself such that the subgroup
\begin{equation}\label{m4.7'.}
\Gamma^{G_t}_z=\big\langle R_{\gamma_1^{G_t}},\ldots ,R_{\gamma_\ell^{G_t}}\big\rangle\subset\mathrm{SU}_{n+1}
\end{equation}
generated by $R_{\gamma_1^{G_t}},\ldots ,R_{\gamma_\ell^{G_t}}$ is dense in $\mathrm{SU}_{n+1}$.
\end{lemma}

\subsection{Perturbed dynamical system}\label{S4.1''}

First, we recall the general Hamiltonian dynamical system $\big(\psi^\mathscr{H}_t,\mathscr{H}^{-1}(c),dv_{\mathscr{H}^{-1}(c)}\big)$ discussed in \eqref{diag1}, and we shall also denote $\psi^\mathscr{H}_t$ by $\psi_t$ for short. When $\mathscr{H}(x,\xi,w)=\Vert\xi\Vert_{T^*X}^2$ and $c=1$, we recover exactly $\big(G_t,q^*S^*X,dv_{q^*S^X}\big)$. Therefore, we assume that $\mathscr{H}$ is a small perturbation of $\Vert\xi\Vert_{T^*X}^2$, then $\mathscr{H}^{-1}(c)\cong q^*S^*X$ through \eqref{1.4..} and we get a perturbed dynamical system.


Similar to \eqref{m4.1.}, we have the following $\psi_t$-invariant decomposition
\begin{equation}
	T\mathscr{H}^{-1}(c)=E^{c,\psi_t}\oplus E^{s,\psi_t}\oplus E^{u,\psi_t}
\end{equation}
and let $W^{s/u,\psi_t}_\mathrm{loc}$ be the corresponding local stable and unstable foliations.

Following Hirsch-Pugh-Shub \cite[Theorems 6.1(f), 6.8]{HPSBook}, we obtain the existance of a perturbed central foliation $\mathcal{F}^{\psi_t}$ similar to \eqref{m4.5..}.
\begin{prop}
	There exists $\tau_0>0$ such that for all $\tau\in[-\tau_0,\tau_0]$, for all $k>0$ sufficiently large and for all $\varepsilon>0$, there exists $\delta>0$ such that if $\psi_t\colon q^*S^*X\to q^*S^*X$ satisfies $d_{C^k(q^*S^*X)}(G_1,\psi_1)\leqslant \delta$, then there exists a central foliation $\mathcal{F}^{\psi_t}$ of $q^*S^*X$ by $C^k$-leaves such that for all $z\in q^*S^*X$,
	\begin{enumerate}
		\item  there exists $J_z\colon \mathcal{F}^{G_t}(z)\to \mathcal{F}^{\psi_t}(z)$ which is a $C^k$-smooth diffeomorphism,
		\item let $I_z\colon\mathcal{F}^{G_t}(z)\to\mathcal{F}^{G_t}(z)$ be the identity map, then $d_{C^k}(J_z,I_z)\leqslant \varepsilon$,
		\item the map $z\mapsto J_z$ is a H\"older continuous map with respect to the $C^k$-norm, where we view $\mathcal{F}^{G_t}(z)=[-\tau,\tau]\times \mathbb{CP}^n$ for $z$ in a small region.
	\end{enumerate}
	In particular, $\mathcal{F}^{\psi_t}$ is a foliation by $C^k$-local leaves whose induced Riemannian volume is bounded uniformly from below and from above.
\end{prop}

We define the perturbed holonomy $H_{zz'}^{s/u,\psi_t}\colon \mathcal{F}^{\psi_t}(z)\to\mathcal{F}^{\psi_t}(z')$ as the unique intersection point
\begin{equation}
H_{zz'}^{s/u,\psi_t}(z'')=W^{s/u}_\mathrm{loc}(z'')\cap \mathcal{F}^{\psi_t}(z'),
\end{equation}
where $z'\in W^{s/u,\psi_t}_\mathrm{loc}$ and $z''\in \mathcal{F}^{\psi_t}(z)$.

Analogous to \eqref{m4.2.}, let $\mathcal{X}_{\mathscr{H}^{-1}}$ be the vector field on $\mathscr{H}^{-1}(c)$ representing the direction of $\psi_t$,
\begin{equation}
	\mathcal{X}_{\mathscr{H}^{-1}(c)}=\frac{d}{dt}\Big|_{t=0}\psi_t.
\end{equation}
Unfortunately, in general, $E^{c,\psi_t}$ \emph{does not} admit a $\psi_t$-invariant decomposition as in \eqref{m4.3.}, nor does it provide a $\psi_t$-invariant $\mathbb{CP}^n$-foliation analogous to \eqref{m4.4.}. This is a major obstacle in our proof. 

The key idea for addressing this absence is that $\mathcal{F}^{\psi_t}(z)$ is tangent to $E^{c,\psi_t}(z)$, which contains $\mathbb{C}\cdot\mathcal{X}_{\mathscr{H}^{-1}(c)}$. Consequently, we can treat $\mathcal{F}^{\psi_t}(z)$ as a saturation by small $\psi_t$-orbit segments
of $J_z(\mathcal{F}(z))$. We then can use the quotient $\big(E^{c,\psi_t}/\mathbb{C}\cdot\mathcal{X}_{\mathscr{H}^{-1}(c)}, \mathcal{F}^{\psi_t}/\psi_t\big)$ to achieve an effect similar to $(\ker q_*,\mathcal{F})$.

\subsection{Reducing to a perturbation of $\mathrm{SU}_{n+1}$-action}\label{S4.2.}


We have the following {\em $k$-bunching} property.
\begin{lemma}\label{claim1}
	For every $k\in\mathbb{N}$ There exists $\delta>0$ such that for all $\psi_t\colon q^*S^*X\to q^*S^*X$ with $d_{C^k(q^*S^*X)}(\psi_1,G_1)\leqslant\delta_0$, we have
	\begin{equation}\label{rBunching}\bV(\psi_1)_*\big|_{E^{s,\psi_t}}\bV\cdot \frac{\bV(\psi_1)_*\big|_{E^{c,\psi_t}}\bV^k}{\bV(\psi_1)_*\big|_{E^{c,\psi_t}}\bV_\mathrm{co}}<1,\end{equation}
	where $\lV\cdot\rV_{\mathrm{co}}$ is the conorm.
\end{lemma}
\begin{pro}
	For $G_t$, we have the strict inequality $\Vert (G_1)_*|_{E^{s,G_t}}\Vert<1$, and $(G_1)_*|_{E^{c,G_t}}$ acts as an isometry, hence \eqref{rBunching} holds for every $k\in\mathbb{N}$. It is easy to see that this property is stable under $C^1$-perturbations.\qed
\end{pro}

Note that hyperbolic splittings can be detected by using cone fields. Carrying a cone field within itself is a $C^1$ open condition, and hence it persists under perturbation. The above $k$-bunching is used to show the following result, see Pugh-Shub-Wilkinson \cite[Page 545]{PSW1} and Hirsch-Pugh-Shub \cite[Theorem 6.7]{HPSBook}. 
\begin{prop}\label{Wilk}
	For any $k>0$ and any $\varepsilon>0$, there exists $\delta>0$, $a>0$, and $C>0$ such that for all flow $\psi_t\colon q^*S^*X\to q^*S^*X$ with $d_{C^k(q^*S^*X)}(\psi_1,G_1)\leqslant \delta$, the holonomy maps $H^{s,\psi_t}_{z_1z_2}\colon\mathcal{F}^{\psi_t}(z_1)\to \mathcal{F}^{\psi_t}(z_2)$ and $H^{s,G_t}_{z_1'z_2'}\colon\mathcal{F}^{G_t}(z_1')\to\mathcal{F}^{G_1}(z_2')$ satisfy
	\begin{equation}\label{shadowing}
		d_{C^k}\Big(J_{z_2}^{-1}H_{z_1z_2}^{s,\psi_t}J_{z_1},H_{z_1'z_2'}^{s,G_1}\Big)\leqslant \varepsilon,
	\end{equation}
	whenever $d(z_1,z_1'),d(z_2,z_2')\leqslant a$, and $d^{s,\psi_t}(z_1,z_2), d^{s,G_t}(z_1',z_2')\leqslant C$, where $d^{s,\psi_t}(\cdot,\cdot),d^{s,G_t}(\cdot,\cdot)$ denotes the distance in the induced Riemannian metric of $W^{s,\psi_t}_\mathrm{loc},W^{s,G_t}_\mathrm{loc}$ respectively. Moreover, a similar statement holds when replacing the roles of $s$ and $u$.
\end{prop}

Now we state a classical result for Anosov geodesic flows.
\begin{lemma}
	There exist $C>0$ and $K\in\mathbb{N}$ such that for all $z\in q^*S^*X$, the loops $\gamma_1,\ldots,\gamma_\ell$ are all composed of at most $K$-many segments, each segment of length at most $C$.
\end{lemma}

The following result shows that we can control the perturbed $su$-path, see \cite[Theorem 6.8]{HPSBook}.
\begin{prop}\label{P4.6}
	For any $\varepsilon>0$, we have $\delta>0$ such that when $d_{C^k(q^*S^*X)}(G_1,\psi_1)\leqslant\delta$, then for every $z\in q^*S^*X$, there exist $\gamma_1^{\psi_t},\ldots,\gamma_\ell^{\psi_t}$ which are $su$-paths from $\mathcal{F}^{\psi_t}(z)$ to itself, all of which are composed of at most $K+\dim S^*X-1$-many segments. In addition,  $\gamma_i^{\psi_t}=(z_1^{\psi_t},\ldots ,z_j^{\psi_t},z'_1,\ldots, z'_{\mathrm{dim}S^*X-1})$ shadows $\gamma_i^{G_t}=(z_1^{G_t},\ldots,z_j^{G_t})$, in the sense that if $z_i^{G_t}$ connects to $z_{i+1}^{G_t}$ through an $s$-path with respect to $G_t$, then $z_i^ {\psi_t} $ connects to $z_{i+1}^ {\psi_t} $ via an $s$-path with respect to $\psi_t$, and $d(z_i^{G_t},z_i^ {\psi_t})\leqslant\varepsilon$. Moreover, the segments which connect $z^{\psi_t}_i$ to $z^{\psi_t}_{i+1}$ are of length at most $(C+1)$, the segments which connect $z_{j}^ {\psi_t} $ to $z'_1$ and $z'_i$ to $z'_{i+1}$ are of length at most $\varepsilon$. A similar statement holds for the path segments where the roles of $s$ and $u$ are reversed.
\end{prop}

For any $su$-path $\gamma^ {\psi_t}=(z_1,\cdots,z_{\vert\gamma^ {\psi_t}\vert})$, similar to \eqref{m4.6.}, we define a $C^k$-smooth map $R_{\gamma^{\psi_t}}\colon\mathcal{F}^{\psi_t}(z_1)\to\mathcal{F}^{\psi_t}\big(z_{\vert\gamma^ {\psi_t}\vert}\big)$ by \begin{equation}\label{m4.13}
	R_{\gamma^{\psi_t}}=H^{s/u,\psi_t}_{z_{\vert\gamma^{\psi_t}\vert-1}z_{\vert\gamma^{\psi_t}\vert}}\cdots H^{s/u,\psi_t}_{z_1z_2}.
\end{equation}
Let us denote
 \begin{equation}\label{m4.14}
\Gamma_z^{\psi_t}=\big\langle R_{\gamma^{\psi_t}_1},\ldots,R_{\gamma^{\psi_t}_{\ell}}\big\rangle,
 \end{equation}
which acts on $\mathcal{F}^{\psi_t}(z)$. By Proposition \ref{P4.6}, $\gamma_i^{\psi_t}$ is composed of at most $(C+\mathrm{dim}S^*X-1)$-many segments, each approximating a segment of $\gamma_i^{G_t}$, together with Proposition \ref{Wilk}, we get
\begin{equation}
	d_{C^k}\big(R_{\gamma^{\psi_t}_i},J_zR_{\gamma_i} J_z^{-1}\big)\leqslant C'\varepsilon.
\end{equation}

Since $H^{s/u,\psi_t}_{zz'}$ commutes with $\psi_t$, we can regard $\Gamma^{\psi_t}_z$ as acting on $\mathcal{F}^{\psi_t}(z)/\psi_t\cong\mathbb{CP}^n$. Therefore, we can view $\Gamma^{\psi_t}_z$ as a perturbation of $\Gamma_z^{G_t}$ given in \eqref{m4.7'.}. Next we shall equip $\mathcal{F}^{\psi_t}(z)/\psi_t$ with a $\Gamma^{\psi_t}_z$-invariant measure. This step is crucially dependent on the Hamiltonian property.

\subsection{Symplectic quotient and leavewise stable ergodicity}

Indeed, we have the following stronger result.
\begin{prop}\label{p4.7}
	There exists symplectic form $\omega_{\mathcal{F}^{\psi_t}(z)/\psi_t}$ on $\mathcal{F}^{\psi_t}(z)/\psi_t$ such that
	\begin{equation}\label{m4.16}
\big(H^{s/u,\psi_t}_{zz'}\big)^*\omega_{\mathcal{F}^{\psi_t}(z')/\psi_t}=\omega_{\mathcal{F}^{\psi_t}(z)/\psi_t}.
	\end{equation}
\end{prop}

\begin{pro}
First, $\omega_{q^*(T^*X)}$ descents to a symplectic form $\omega_{\mathscr{H}^{-1}(c)/\psi_t}$ locally on $\mathscr{H}^{-1}(c)/\psi_t$. Indeed, this is a toy version of Marsden-Weinstein symplectic quotient, see for instance McDuff-Salamon \cite[\S\,5.4]{MR3674984}. It relies on the definition of the Hamiltonian flow,
\begin{equation}\label{m4.17}
	\iota_{\mathcal{X}_{\mathscr{H}^{-1}(c)}}\big(i_c^*\omega_{q^*(T^*X)}\big)=d\mathscr{H}=0,
\end{equation}
where $i_c\colon \mathscr{H}^{-1}(c)\to q^*(T^*X)$ is the natural embedding used in \eqref{diag1} and the second equality holds because $\mathscr{H}$ is a constant on $\mathscr{H}^{-1}(c)$.

We can give a more concrete description. Let us take a representative $S_z\subset\mathcal{F}^{\psi_t}(z)$ of $\mathcal{F}^{\psi_t}(z)/\psi_t$, for example $J_z(\mathcal{F}(z))$, then we pullback $i_c^*\omega_{q^*(T^*X)}$ throught $S_z\hookrightarrow \mathscr{H}^{-1}(c)$ to get $\omega_{S_z}$. Any other representative $S_{z}'$ can be written as the image of a map
\begin{equation}
f(z_1)=\psi_{t(z_1)}(z_1),
\end{equation}
where $z_1\in S_z$ and $t(\cdot)$ is a function on $S_z$. We have $f^*\omega_{S_z'}=\omega_{S_z}$. To see this, first note that $\omega_{q^*T^*X}$ is $\psi_t$-invariant, hence at the point $f(z_1)$, we have $\psi_{t(z_1)}^*\omega_{S_z'}=\omega_{S_z}$, and the difference $f^*\omega_{S_z'}-\psi_{t(z_1)}^*\omega_{S_z'}$ involves the the contraction of $i_c^*\omega_{q^*(T^*X)}$ with the direction of $\psi_t$, which vanishes by \eqref{m4.17}.

Now we turn to \eqref{m4.16}. As the local leaf of the central foliation, $\mathcal{F}^{\psi_t}$ is $\psi_t$-invariant in the sense that $\mathcal{F}^{\psi_t}(\psi_s(z))=\psi_s(\mathcal{F}^{\psi_t}(z))$. Consequently, if two points in $\mathcal{F}^{\psi_t}(z)$ and in $\mathcal{F}^{\psi_t}(z')$ which are connected via a stable or unstable leaf, we may push them forward by $\psi_t$, under which $\omega_{\mathcal{F}^{\psi_t}(z)/\psi_t}$ and $\omega_{\mathcal{F}^{\psi_t}(z')/\psi_t}$ remains invariant, while the two points become arbitrarily close. \qed
\end{pro}

By \eqref{m4.13}, \eqref{m4.14} and \eqref{m4.16}, we obtain the following result.
\begin{coro}
The $\Gamma_z^{\psi_t}$-action preserves $\omega_{\mathcal{F}^{\psi_t}(z)/\psi_t}$,  thus also the induced volume $dv_{\mathcal{F}^{\psi_t}(z)/\psi_t}=\omega_{\mathcal{F}^{\psi_t}(z)/\psi_t}^n$.
\end{coro}

\begin{remark}
By the construction in the proof of Proposition \ref{p4.7}, we form a volume $dv_{\mathcal{F}^{\psi_t}(z)}=\omega_{S_z}\wedge dt$ on $\mathcal{F}^{\psi_t}(z)$. It is independent of the representative $S_z$ and is invariant with respect to both $\psi_t$ and $\Gamma_z^{\psi_t}$.
\end{remark}

Now since $dv_{\mathcal{F}^{\psi_t}(z)/\psi_t}$ is induced by the symplectic form $\omega_{\mathcal{F}^{\psi_t}(z)/\psi_t}$, by applying the Moser's trick, see for instance \cite[Lemma 3.2.1]{MR3674984}, we can assume that $dv_{\mathcal{F}^{\psi_t}(z)/\psi_t}$ is a constant multiple of the Fubini-Study volume on $\mathbb{CP}^n$ without loss of generality. Under these conditions the group action stable ergodicity results of Dolgopyat-Krikorian \cite[Corollary 2]{DimaKrikorian} and DeWitt \cite[Theorem 1]{MR4756948} apply, yielding the following leavewise stable ergodicity.
	\begin{prop}\label{Dol}
		There exist $k\in\mathbb{N}$ and $\delta>0$ such that if $d_{C^k(q^*S^*X)}(G_1,\psi_1)\leqslant \delta$, then for all $z\in q^*S^*X$, the $\Gamma_z^{\psi_t}$-action $(\mathcal{F}^{\psi_t}(z)/\psi_t,dv_{\mathcal{F}^{\psi_t}(z)/\psi_t})$ is ergodic.
	\end{prop}

\subsection{Flow stable ergodicity}\label{S4.3.}

First, we state the following important result of Pugh-Shub \cite{Juliennes}.
\begin{prop}\label{Juli}
	For any $k>0$ and any volume preserving $\psi_t\colon q^*S^*X\to q^*S^*X$ such that $d_{C^k(q^*S^*X)}(\psi_1,G_1)$ is sufficiently small (hence with bunching), for every $\psi_1$-invariant set $B$ such that 
	\begin{enumerate}
		\item $B$ is $s$-saturated modulo the volume,
		\item $B$ is $u$-saturated modulo the volume,
	\end{enumerate}  
	$B$ must be $su$-saturated modulo the volume, where $su$-saturation means saturation by $su$-paths.
\end{prop}

We have the following classical Anosov absolute continuity.
\begin{prop}\label{Anosov}
	For any $z\in q^*S^*X$, let $B$ be a measurable subset of $\mathcal{F}^{\psi_t}(z)$, then $dv_{\mathcal{F}^{\psi_t}(z)}(B)=0$ if and only if $dv_{\mathscr{H}^{-1}(c)}(\mathrm{Acc}(B))=0$, where $\mathrm{Acc}(B)$ is the accessibility class of $B$, that is, the saturation by $su$-paths.
\end{prop}

Now we are ready to prove our main dynamical result, the ergodicity of $\psi_t$.
\begin{theo}\label{T6.3}
	There exist $k\in\mathbb{N}$ and $\delta>0$, such that for $d_{C^k(q^*S^*X)}(G_1,\psi_1)\leqslant \delta$, $(\psi_t,dv_{\mathscr{H}^{-1}(c)})$ is ergodic. In particular, $(G_t,dv_{q^*S^*X})$ is ergodic.
\end{theo}
\begin{pro}
	For $\mathscr{A}\in C(q^*S^*X)$, we set
	\begin{equation}
		\langle\mathscr{A}\rangle_1=\int_{0}^1\mathscr{A}\circ\psi_tdt.
	\end{equation}
	Since $\psi_t$ preserves $dv_{\mathscr{H}^{-1}(c)}$, the limits
	\begin{equation}
		\lim_{j\to \infty} \frac{1}{j}\sum_{0\leqslant i\leqslant j}\langle\mathscr{A}\rangle_1\circ \psi_1^i,\ \ \lim_{j\to \infty} \frac{1}{j}\sum_{-j\leqslant i\leqslant 0}\langle\mathscr{A}\rangle_1\circ \psi_1^i
	\end{equation}
	are well-defined for $dv_{\mathscr{H}^{-1}(c)}$-a.e.
	
	For any $a\in \mathbb{R}$, let us define
	\begin{equation}
		\begin{split}
			&B^s=\Big[\lim_{j\to \infty} \frac{1}{j}\sum_{0\leqslant i\leqslant j}\langle\mathscr{A}\rangle_1\circ \psi_1^i\leqslant a\Big],\ \ B^u=\Big[\lim_{j\to \infty} \frac{1}{j}\sum_{-j\leqslant i\leqslant 0}\langle\mathscr{A}\rangle_1\circ \psi_1^i\leqslant a\Big],\\
			&B=B^s\cap B^u.
		\end{split}
	\end{equation}
	By the pointwise ergodic theorem for every ergodic component of $dv_{\mathscr{H}^{-1}(c)}$, it is clear that $B^s$ is $s$-saturated, $B^u$ is $u$-saturated, and both are equal to $B$ modulo the volume. Hence $B$ satisfies the assumptions of Theorem \ref{Juli}, and hence there exists $B'$ which is $su$-saturated and $dv_{\mathscr{H}^{-1}(c)}(B\Delta B')=0$. Moreover, we may assume that $B'$ is $\psi_t$-invariant, as $B$ is $\psi_t$-invariant.
	
	By Theorem \ref{Anosov}, $dv_{\mathscr{H}^{-1}(c)}(B')>0$ if and only if $dv_{\mathcal{F}^{\psi_t}(z)}(B'\cap\mathcal{F}^{\psi_t})>0$ for all $z$, otherwise if the intersection is of $0$ induced volume, since $B'$ is its own accessibility class, $B'$ must be of $0$ volume. We rely on the observation that $B'$ must intersect in a non-empty way every $\mathcal{F}^{\psi_t}(z)$, $z\in q^*S^*X$. This follows from the following observation: $g_1$ is an accessible map on $S^*X$, which admits an $su$-path between any two points in $S^*X$, composed of at most $K$-many segments, each of length at most $C$. Perturbing $G_1$ a bit induces a new $su$-path for $\psi_1$, which is between a small neighborhood of any $z$, and a small neighborhood of any $\mathcal{F}^{\psi_t}(z')$, for any $z,z'\in q^*S^*X$. 
		Therefore, 
	$B'\cap \mathcal{F}^{\psi_t}(z)\neq \varnothing$ for every $z\in q^*S^*X$. 
	
	Assume that $dv_{\mathscr{H}^{-1}(c)}(B')>0$. Let $z\in q^*S^*X$, and note that $B'\cap \mathcal{F}^{\psi_t}(z)$ is invariant under holonomies into $\mathcal{F}^{\psi_t}(z)$, hence it is $\Gamma_z^{\psi_t}$-invariant. Moreover, since $B'$ is $\psi_t$-invariant, and $\psi_t$ preserves $dv_{\mathscr{H}^{-1}(c)}$, by Proposition \ref{Dol}, $B'$ has a full measure in $\mathcal{F}^{\psi_t}(z)$, hence its complement has a zero measure, hence $(B')^c$ has zero volume.
	.
	
	It follows that $dv_{\mathscr{H}^{-1}(c)}(B)\in \{0,1\}$, for every $a\in \mathbb{R}$ and every $\mathscr{A}\in C(q^*S^*X)$. It is an easy consequence of the pointwise ergodic theorem that in that case $(\psi_t,dv_{\mathscr{H}^{-1}(c)})$ is ergodic.\qed
\end{pro}

\section{Moment map and applications}\label{S5..}

In this section, we present more general examples to which Theorems \ref{mt1.1}, \ref{mt1.2} and \ref{mt1.4} apply.

	\subsection{Moment map}\label{s6.1}
	
	We use the notation as in \cref{s1.2} and \cref{Ca}. In particular, $L$ is a holomorphic line bundle over the complex manifold $N$.

	Let $U$ be a compact connected Lie group with its Lie algebra $\mathfrak{u}$. We assume that $U$ acts holomorphically on $N$, and this action lifts to a holomorphic unitary action on $L$. For $\mathfrak{a}\in\mathfrak{u}$, let $\mathcal{L}^L_\mathfrak{a}$ denote the Lie derivative of $\mathfrak{a}$ on ${C}^\infty(N,L)$, that is,
	\begin{equation}\label{m4.1}
		\mathcal{L}^{L}_\mathfrak{a}s(w)=\frac{d}{dt}\Big|_{t=0}e^{t\mathfrak{a}}\big(s(e^{-t\mathfrak{a}}w)\big).
	\end{equation}
	Let $\nabla^L$ be the Chern connection on $L$, then the Kostant formula \cite[Definition 7.5]{MR2273508} gives a moment map
	$\mu_L\colon N\rightarrow\mathfrak{u}^*$ by
	\begin{equation}\label{hc1}
		{2\pi i}\langle\mu_L,\mathfrak{a}\rangle=\nabla^L_{\mathfrak{a}}-\mathcal{L}^L_{\mathfrak{a}}.
	\end{equation}
	
	Since $U$ acts holomorphically and unitarily on $\det T^{(1,0)}N$, the determinant line bundle of the holomorphic tangent bundle $T^{(1,0)}N$ on $N$, we also have an associated moment map and denote it by $\mu_{\det T^{(1,0)}N}\colon N\rightarrow \mathfrak{u}^*$.

	Similar to \eqref{m4.1}, let $\mathcal{L}^{L^p}_\mathfrak{a}$ denote the Lie derivative on ${C}^\infty(N,L^p)$, then Bismut-Ma-Zhang \cite[Theorem 3.1]{MR3615411} asserts that when restricts to the holomorphic sections $H^{(0,0)}(N,L^p)$, $p^{-1}\mathcal{L}^{L^p}_\mathfrak{a}$ is a Toeplitz operator in the sense of \eqref{cb16}.
	\begin{prop}
		For any $\mathfrak{a}\in\mathfrak{u}$, we have
		\begin{equation}\label{hc6}
			p^{-1}\mathcal{L}_\mathfrak{a}^{L^p}\big|_{H^{(0,0)}(N,L^p)}=-{2\pi \sqrt{-1}} T_{\langle\mu_L+p^{-1}\mu_{\det T^{(1,0)}N},\mathfrak{a}\rangle,p}.
		\end{equation}
	\end{prop}

	\subsection{More examples}\label{Sa4.2}
	
	Now we suppose further that the $\pi_1(X)$-action in \eqref{m1.1} is induced by a representation
	\begin{equation}\label{m5.4}
		\rho\colon\pi_1(X)\longrightarrow U.
	\end{equation}

	Let $U(\mathfrak{u})$ be the universal enveloping algebra of $\mathfrak{u}$, and $U(\mathfrak{u})^{\leqslant j}\subset U(\mathfrak{u})$ the subspace generated by monomials with degree less or equal to $j$. Using \eqref{hc6}, we can define a principal symbol map
	\begin{equation}\label{m4.5.}
		\begin{split}
			\sigma_j\colon U(\mathfrak{u})^{\leqslant j}&\to{C}^\infty(N)\\
			\mathfrak{a}_1\cdots\mathfrak{a}_j&\to (-1)^j\big(2\pi\sqrt{-1}\langle\mu_L,\mathfrak{a}_1\rangle\big)\cdots \big(2\pi\sqrt{-1}\langle\mu_L,\mathfrak{a}_j\rangle\big).
		\end{split}
	\end{equation}
	We set
	\begin{equation}
		\mathscr{U}(\mathfrak{u})^{\leqslant j}=\pi_1(X)\backslash\big(\widetilde{X}\times U(\mathfrak{u})^{\leqslant j}\big),
	\end{equation}
	then \eqref{m4.5.} extends to a map
	\begin{equation}
		\sigma_j\colon \mathscr{U}(\mathfrak{u})^{\leqslant j}\to {C}^\infty(\mathscr{N})
	\end{equation}
	where $\mathscr{N}$ is defined in \eqref{diag1}.
	
	Let us choose sequences $\big\{\Theta_j\in {C}^\infty\big(X,T^*X\otimes \mathscr{U}(\mathfrak{u})^{\leqslant j}\big)\big\}_{j=0}^\ell$ and $\big\{V_j\in{C}^\infty\big(X,\mathscr{U}(\mathfrak{u})^{\leqslant j}\big)\big\}_{j=0}^k$, then we define
	\begin{equation}
		\Theta^{F_p}=\sum_{j=0}^{\ell}p^{-j}\Theta_j,\ \  V^{F_p}=\sum_{j=0}^{k}p^{-j}V_j.
	\end{equation}
	Since the action of $\mathfrak{u}$ on $H^{(0,0)}(N,L^p)$ induces an action of $\mathscr{U}(\mathfrak{u})^{\leqslant j}$ on $F_p$, we can view 
	\begin{equation}
		\Theta^{F_p}\in {C}^\infty(X,T^*X\otimes \mathrm{End}(F_p)),\ \ V^{F_p}\in {C}^\infty(X,\mathrm{End}(F_p)).
	\end{equation}
	Let $\{e_i\}_{i=1}^r$ be a local orthonormal frame of 
	$TX$, then we define a perturbed Hamiltonian
	\begin{equation}\label{m4.7}
		\begin{split}
			\mathscr{H}^{F_p}=&\sum_{i=1}^{r}-\big(p^{-1}\nabla^{F_p}_{e_i}-\Theta^{F_p,*}(e_i)\big)\big(p^{-1}\nabla^{F_p}_{e_i}+\Theta^{F_p}(e_i)\big)\\
			&+p^{-1}\big(p^{-1}\nabla^{F_p}_{\nabla^{TX}_{e_i}e_i}+\Theta^{F_p}(\nabla^{TX}_{e_i}e_i)\big)+V^{F_p}+V^{F_p,*}.
		\end{split}
	\end{equation}
	where $\Theta^{F_p}$ plays the role of a field and $V^{F_p}$ a potential.

	We define the principal symbol $\Theta\in{C}^\infty(\mathscr{N},q^*TX)$ of $\Theta^{F_p}$ and the principal symbol $V\in{C}^\infty(\mathscr{N})$ of $V^{F_p}$ by
	\begin{equation}
		\Theta=\sum_{i=1}^{\ell}\sigma_i(\Theta_i),\ \ \  V=\sum_{i=1}^{k}\sigma_i(V_i).
	\end{equation}
	The principal symbol $\mathscr{H}$ of $\mathscr{H}^{F_p}$ is given by
	\begin{equation}
		\mathscr{H}=-g^{ij}\big(\sqrt{-1}\xi_i-\overline{\Theta(\pa_{x_i})}\big)\big(\sqrt{-1}\xi_j+\overline{\Theta(\pa_{x_j})}\big)+V+\overline{V}.
	\end{equation}
	Hence \eqref{m3.14;} is verified, and Theorem \ref{mt1.1} can be applied to $\mathscr{H}^{F_p}$ given in \eqref{m4.7}.



	
Now we restrict to the case
\begin{equation}\label{}
	\big(X,N,L,H^{(0,0)}(N,L^p),U,\mathfrak{u}\big)=\big(\Gamma\backslash\mathbb{H}^r,\mathbb{CP}^n,O_{\mathbb{CP}^n}(1),\mathrm{Sym}^p\mathbb{C}^{n+1},\mathrm{SU}_{n+1},\mathfrak{su}_{n+1}\big).
\end{equation}
We refer to Ma-Ma \cite[\S\,5.1]{MR4808253} for  examples of $\rho\colon \Gamma\to\mathrm{SU}_{n+1}$ with dense image, both generic and arithmetic. If $\lv\Theta_i\rv_{C^k(X)},\lv V_i\rv_{C^k(X)}$ are small, then Theorem \ref{mt1.4} is applicable. In particular, this generalizes the example in \eqref{1.22}.

\addcontentsline{toc}{section}{References}
\bibliographystyle{abbrv}

\begin{thebibliography}{10}
	
	\bibitem{MR2273508}
	N.~Berline, E.~Getzler, and M.~Vergne.
	\newblock {\em Heat kernels and {D}irac operators}.
	\newblock Grundlehren Text Editions. Springer-Verlag, Berlin, 2004.
	\newblock Corrected reprint of the 1992 original.
	
	\bibitem{MR2838248}
	J.-M. Bismut, X.~Ma, and W.~Zhang.
	\newblock Op\'{e}rateurs de {T}oeplitz et torsion analytique asymptotique.
	\newblock {\em C. R. Math. Acad. Sci. Paris}, 349(17-18):977--981, 2011.
	
	\bibitem{MR3615411}
	J.-M. Bismut, X.~Ma, and W.~Zhang.
	\newblock Asymptotic torsion and {T}oeplitz operators.
	\newblock {\em J. Inst. Math. Jussieu}, 16(2):223--349, 2017.
	
	\bibitem{MR658304}
	R.~Bott and L.~W. Tu.
	\newblock {\em Differential forms in algebraic topology}, volume~82 of {\em
		Graduate Texts in Mathematics}.
	\newblock Springer-Verlag, New York-Berlin, 1982.
	
	\bibitem{MR1717580}
	K.~Burns and A.~Wilkinson.
	\newblock Stable ergodicity of skew products.
	\newblock {\em Ann. Sci. \'Ecole Norm. Sup. (4)}, 32(6):859--889, 1999.
	
	\bibitem{cekić2024semiclassicalanalysisprincipalbundles}
	M.~Cekić and T.~Lefeuvre.
	\newblock Semiclassical analysis on principal bundles.
	\newblock {\em arXiv: 2405.14846}, 2024.
	
	\bibitem{MR818831}
	Y.~Colin~de Verdi{\`e}re.
	\newblock Ergodicit\'{e} et fonctions propres du laplacien.
	\newblock {\em Comm. Math. Phys.}, 102(3):497--502, 1985.
	
	\bibitem{MR4756948}
	J.~DeWitt.
	\newblock Simultaneous linearization of diffeomorphisms of isotropic manifolds.
	\newblock {\em J. Eur. Math. Soc. (JEMS)}, 26(8):2897--2969, 2024.
	
	\bibitem{DimaKrikorian}
	D.~Dolgopyat and R.~Krikorian.
	\newblock On simultaneous linearization of diffeomorphisms of the sphere.
	\newblock {\em Duke Math. J.}, 136(3):475--505, 2007.
	
	\bibitem{MR3969938}
	S.~Dyatlov and M.~Zworski.
	\newblock {\em Mathematical theory of scattering resonances}, volume 200 of
	{\em Graduate Studies in Mathematics}.
	\newblock American Mathematical Society, Providence, RI, 2019.
	
	\bibitem{HPSBook}
	M.~W. Hirsch, C.~C. Pugh, and M.~Shub.
	\newblock {\em Invariant manifolds}, volume Vol. 583 of {\em Lecture Notes in
		Mathematics}.
	\newblock Springer-Verlag, Berlin-New York, 1977.
	
	\bibitem{MR2093043}
	D.~Huybrechts.
	\newblock {\em Complex geometry}.
	\newblock Universitext. Springer-Verlag, Berlin, 2005.
	\newblock An introduction.
	
	\bibitem{Ma-Ma}
	M.~Ma and Q.~Ma.
	\newblock Semiclassical analysis, geometric representation and quantum
	ergodicity.
	\newblock {\em Communications in Mathematical Physics}, 405(11):259, 2024.
	
	\bibitem{MR4808253}
	M.~Ma and Q.~Ma.
	\newblock Semiclassical analysis, geometric representation and quantum
	ergodicity.
	\newblock {\em Comm. Math. Phys.}, 405(11):Paper No. 259, 28, 2024.
	
	\bibitem{MR4665497}
	Q.~Ma.
	\newblock Toeplitz operators and the full asymptotic torsion forms.
	\newblock {\em J. Funct. Anal.}, 286(3):Paper No. 110210, 74pp, 2024.
	
	\bibitem{MR2339952}
	X.~Ma and G.~Marinescu.
	\newblock {\em Holomorphic {M}orse inequalities and {B}ergman kernels}, volume
	254 of {\em Progress in Mathematics}.
	\newblock Birkh\"{a}user Verlag, Basel, 2007.
	
	\bibitem{MR2393271}
	X.~Ma and G.~Marinescu.
	\newblock Toeplitz operators on symplectic manifolds.
	\newblock {\em J. Geom. Anal.}, 18(2):565--611, 2008.
	
	\bibitem{MR3674984}
	D.~McDuff and D.~Salamon.
	\newblock {\em Introduction to symplectic topology}.
	\newblock Oxford Graduate Texts in Mathematics. Oxford University Press,
	Oxford, third edition, 2017.
	
	\bibitem{MR4611826}
	M.~Puchol.
	\newblock The asymptotics of the holomorphic analytic torsion forms.
	\newblock {\em J. Lond. Math. Soc. (2)}, 108(1):80--140, 2023.
	
	\bibitem{Juliennes}
	C.~Pugh and M.~Shub.
	\newblock Stable ergodicity and julienne quasi-conformality.
	\newblock {\em J. Eur. Math. Soc. (JEMS)}, 2(1):1--52, 2000.
	
	\bibitem{PSW1}
	C.~Pugh, M.~Shub, and A.~Wilkinson.
	\newblock H\"older foliations.
	\newblock {\em Duke Math. J.}, 86(3):517--546, 1997.
	
	\bibitem{MR0995750}
	R.~Schrader and M.~E. Taylor.
	\newblock Semiclassical asymptotics, gauge fields, and quantum chaos.
	\newblock {\em J. Funct. Anal.}, 83(2):258--316, 1989.
	
	\bibitem{MR0402834}
	A.~I. \v{S}nirel$'$man.
	\newblock Ergodic properties of eigenfunctions.
	\newblock {\em Uspehi Mat. Nauk}, 29(6(180)):181--182, 1974.
	
	\bibitem{MR916129}
	S.~Zelditch.
	\newblock Uniform distribution of eigenfunctions on compact hyperbolic
	surfaces.
	\newblock {\em Duke Math. J.}, 55(4):919--941, 1987.
	
	\bibitem{MR1183602}
	S.~Zelditch.
	\newblock On a ``quantum chaos'' theorem of {R}. {S}chrader and {M}. {T}aylor.
	\newblock {\em J. Funct. Anal.}, 109(1):1--21, 1992.
	
	\bibitem{MR2952218}
	M.~Zworski.
	\newblock {\em Semiclassical analysis}, volume 138 of {\em Graduate Studies in
		Mathematics}.
	\newblock American Mathematical Society, Providence, RI, 2012.
	
\end{thebibliography}

\def\cprime{$'$} \def\cprime{$'$}

\end{document}